\documentclass[12pt]{article}
\usepackage{amsmath,amsthm, amssymb }
\newtheorem{theorem}{\noindent Theorem}
\newtheorem{lemma}{\noindent Lemma}[section]
\newtheorem{corollary}{\noindent Corollary}
\newtheorem{proposition}{\noindent Proposition}

\theoremstyle{remark}
\newtheorem{remark}{\it Remark}[section]

\makeatletter

\@addtoreset{equation}{section}
\makeatother

\def\Z{\mathbb{Z}}

\newcommand{\ep}{\varepsilon} 
\newcommand{\de}{\delta}

\newcommand{\la}{\lambda}
\newcommand{\fa}{\varphi}

\newcommand{\La}{\Lambda}

\newcommand{\beq}{\begin{eqnarray*}}
\newcommand{\eeq}{\end{eqnarray*}}
\newcommand{\beqn}{\begin{equation}}
\newcommand{\eeqn}{\end{equation}}

\begin{document}

\begin{center}
{\Large   The potential function and  ladder heights  \\
of a recurrent  random walk on $\Z$
with infinite variance}
\vskip6mm
{K\^ohei UCHIYAMA} \\
\vskip2mm
{Department of Mathematics, Tokyo Institute of Technology} \\
{Oh-okayama, Meguro Tokyo 152-8551\\
}
\end{center}

\vskip6mm

\begin{abstract}
We consider  a recurrent  random walk of i.i.d. increments on the one-dimensional integer lattice
and obtain a formula relating  the hitting distribution of a half-line with the potential function, $a(x)$, of the random walk.  Applying it, we derive an asymptotic estimate of $a(x)$ and thereby a criterion for $a(x)$ to be bounded on a half-line. The application is also made to estimate some hitting probabilities as well as to derive asymptotic behaviour  for large times of the walk conditioned never  to visit the origin.
 \footnote{
{\it key words}: recurrent random walk;  ladder height; potential function; infinite variance; first hitting time }
\end{abstract}

\vskip6mm

\section{Introduction}

Let $S_n=S_0+ X_1+\cdots+ X_n$ be a  random walk on  $\mathbb{Z}$  where the starting position  $S_0$ is an unspecified integer and  the increments $ X_1,  X_2, \ldots$ are  independent and identically distributed  random variables  defined on some probability  space $(\Omega, {\cal F}, P)$ and taking values in $\mathbb{Z}$.  Let $ X$ be a random variable having the same law as $ X_1$.  We suppose throughout the paper that

\vskip2mm
\quad  the walk $S_n$ is \textit{recurrent} and {\it irreducible} (as a Markov chain on $\mathbb{Z}$). 

\vskip2mm\noindent
For a subset $B$ of the whole real line $\mathbb{R}$ such that $B\cap \mathbb{Z} \neq \emptyset$,   put  $\sigma_{B}=\inf \{n\geq 1: S_n\in B\}$,  the first entrance time of  the walk into $B$. Let $Z$ be the first strictly ascending ladder height that is  defined   by
 \[
 Z = S_{\sigma_{[S_0+1,\infty)}}-S_0.
 \]
  We also define  $\hat Z =  S_{\sigma_{(-\infty, S_0-1]}} - S_0$,  the first strictly descending ladder height.
Because of recurrence of the walk  $Z$ is a proper random variable whose distribution is  concentrated on positive integers $x=1, 2, \ldots$ and similarly for  $-\hat Z$.   Let $E$  indicate the integration   by $P$  as usual.
If $\sigma^2:=E  X^{2}<\infty$, then $E Z<\infty$, whereas  if $\sigma^2=\infty$, either $E Z=\infty $ or $E |\hat  Z|=-\infty $ (cf \cite[Section 17]{S}, \cite[Theorem 8.4.7]{Ch}). 

Denote by  $P_x$  the probability of the  random walk  with $S_0=x$ and $E_x$ the expectation by $P_x$.  Put $p^n(x)=P_0[S_n=x]$, $p(x)=p^1(x)$ and define 
 \[
 a(x)=\sum_{n=0}^\infty\big[p^n(0)-p^n(-x)\big];
 \]
 the series on the RHS is convergent  (cf. Spitzer \cite[P28.8]{S}).   The function $a(x)$, called  potential function,  plays a central role in the potential theory of recurrent random walks.
   (This is true for two dimensional walks but here we restrict our discussion below to the one dimensional walks).  
  Spitzer \cite{S0} established  fundamental facts concerning  $a(x)$---its existence, positivity, asymptotic behaviour etc.---and based on them  Kesten and Spitzer  \cite{KS}   obtained certain ratio limit theorems  for the distributions of  the hitting times and sojourn times  of a finite set and the transition probabilities of the walk stopped as it hits the set, which were refined by Kesten \cite{K_RLT}  under  mild additional assumptions.  
 An excellent exposition for the principal contents of   \cite{S0}, \cite{KS}, \cite{K_RLT} is  given in Chapter 7 of Spitzer's book \cite{S};  extensions to non-lattice random walks are obtained by Ornstein \cite{Or},  
 Port and Stone \cite{PS},  Stone \cite{St}. Kesten \cite{K3} conjectured that  the series that defines $a(x)$  converges absolutely and 
 provided certain  mild sufficient conditions for the absolute convergence.
 
 According to  Theorems 6a and  7 of \cite{K_RLT}  
  \beqn\label{RLT}
 \lim_{n\to\infty} \frac{P_x[ S_n=y, S_k\neq 0 \;\;  \mbox{for}\;\;  1\leq k< n]}{P_0[S_n=0,  S_k\neq 0 \;\; \mbox{for}\, 1\leq k< n]} = a^\dagger(x)a^\dagger(-y)+\frac{xy}{\sigma^4}  
  \eeqn
($x,y \in \mathbb{Z}$).
  Here (and in the sequel) $a^\dagger(0)=1$ and $=a(x)$ for $x\neq 0$ and  $1/\infty$ is understood to be zero.
 The  asymptotic estimates  valid uniformly in $x$ and $y$
 of the ratio under the limit  above  are studied  by the present author in \cite{U1dm}, \cite{U1dm.f.s} in case     $\sigma^2 <\infty$    
 and  in \cite{Uattrc} for the stable walks with exponent   $1<\alpha <2$.  
 The denominator  of the ratio in (\ref{RLT}), which  equals  the probability  that the walk starting at zero returns to zero at  $n$  for the first time,   are estimated with some exact asymptotics   in these articles. 
   
  The  following basic  properties of  $a(x)$ are found  in \cite{S}:
\beqn\label{a-a}
a(x +1) -a(x) \to \pm 1/\sigma^2 \quad\mbox{as}\quad x\to \pm \infty
\eeqn
 and 
  \beqn\label{Int_1}
 a(x)-\frac{x}{\sigma^2} \left\{ \begin{array}{ll} =0 \;\; \mbox{for all} \;\; x>0 \quad \mbox{if} \quad P[X\leq -2] =0, \\
 >0 \;\; \mbox{for all} \;\; x>0 \quad \mbox{otherwise}
 \end{array}\right.
 \eeqn 
 (the strict positivity  in the second case of (\ref{Int_1})  is implicit in  \cite{S} if $\sigma^2<\infty$; see e.g., \cite[Eq(2.9)]{U1dm}). When $\sigma^2 <\infty$ (\ref{a-a}) entails the exact asymptotics $a(x)\sim |x|/\sigma^2$, whereas  in case  $\sigma^2 = \infty$ it gives only $a(x)=o(|x|)$ and sharper asymptotic estimates are desired. For the stable walks exact results are given in \cite{B2} for the case $1<\alpha <2$ apart from an extreme case and in \cite[Section 8.1.1]{Upot} 
 for all the cases  $1\leq \alpha \leq 2$ under some natural side conditions.  
 
  In the general case of $\sigma^2 =\infty$ there seems to have been no results of asymptotic estimates of $a(x)$ other than those mentioned above. Very recently the present author gave some relevant results.
Let  $\sigma^2 =\infty$ and  $E|X|<\infty$ and put
\beqn\label{fa}
 m_-(x)={\int_0^x dy \int_y^\infty P[ X<-u]du},\quad  m_+(x)={\int_0^x dy \int_y^\infty P[ X>u]du}
\eeqn
and $m(x) = m_-(x) + m_+(x)$.  It is shown  in \cite{Upot} that   $a(x)+a(-x) \geq C_* x/m(x)$ ($x>0$) with a universal  constant $C_*>0$;   the upper bound is also given so that 
   \beqn\label{a/m}
   a(x) +a(-x)\asymp x/m(x)
   \eeqn
      under a reasonable  side condition, which is satisfied if e.g.,  
 $\limsup_{x\to\infty}  xm'(x)/m(x) <1$ or $m_+(x)/m(x) \to 0$. (Here  $b_x\asymp c_x$ means that $b_x/c_x$ is bounded away from zero and infinity.) 

 In this paper we shall show, supposing $\sigma^2=\infty$, that 
   \beqn\label{Int_4}
\left\{ \begin{array}{ll}   \; {\displaystyle a(x)/ V_{\rm ds}(x) \to 1/EZ \;\;\mbox{and}\;\;  
 E_x[a(S_{\sigma(-\infty,0]})]/a(x) \to 0 }\;\;\ (x\to \infty) \quad  &\mbox{if}  \quad  EZ<\infty,\\[2mm]
{\displaystyle \liminf_{x\to +\infty}  a(x)/V_{\rm ds}(x) =0} \;\; \mbox{and}\;\;
 a(x)= E_x[a(S_{\sigma(-\infty,0]})] \;(x>0) \;
  \quad &\mbox{  otherwise,}
  \end{array}\right.
  \eeqn
  where $V_{\rm ds}$ denotes the renewal function for the weakly descending ladder process, and that  there exists $\lim_{x\to\, \infty} a(-x) \leq \infty$, and 
  \beqn\label{Int_3}
  \;   \begin{array}{ll}
 0< {\displaystyle \lim_{x\to \, \infty} a(-x) <\infty }\quad 
  &\mbox{if} \quad \left\{ \begin{array}{ll} E|X|< \infty \quad \mbox{and}\\
  \int_0^\infty [t/m_-(t)]^2P[X>t] dt <\infty,
  \end{array}\right.\\[4mm]
{\displaystyle \liminf_{x\to \, \infty} a(-x)  = \infty}  \quad &\mbox{otherwise,}
 \end{array}
 \eeqn
 provided  $P[X\geq 2]>0$. 
In (\ref{Int_4}) the $\liminf$ may be expected to be  replaced by $\lim$ (see Remark \ref{rem2.1}(e)). 
Note that if $E|X| = \infty $,  then  $EX_+ = E X_- = \infty$  because of the assumed recurrence of the walk. (Here $X_+=\max\{X, 0\}$ and $X_- = X_+ -X$.)    Applying  (\ref{Int_4}),  we derive  asymptotic estimates of some hitting probabilities as well as  asymptotic behaviour  for large times of the random walk conditioned never  to visit the origin. As an intelligible  manifestation of the significance of the condition $E Z<\infty$ in the sample path behavior of the walk,  we shall observe  that   $E Z<\infty$ if and only if   the walk conditioned never to visit the origin  approaches  the positive infinity with probability one  (Section 7).

 The main results (Theorems \ref{thm1.1} and \ref{thm1.2}) of the present paper are derived from those given in Spitzer's book (that are stated in Section 3 of the present paper)  independently of those of  \cite{Upot} in which the proof is solely based on the Fourier integral representation of $a(x)$.
In the proof of our main results, we could apply those from \cite{Upot} whose usage, however,  we avoid in order  not to cause any  suspicion of circular arguments, some of our results (\ref{Int_4}) being used in \cite{Upot}.
. For the sake of comparison, we  include  the case of finite variance when   all the results are  known or  easily derived from known ones.
\section{Statements of results}
Let $S_n$ be the random walk specified in Introduction and     $Z, \hat Z,  \sigma_B, m_\pm(x)$ and $a(x)$ be as given there.
In order to state the results of the paper we further bring in the following  notation. 
Put 
\[
T=\sigma_{(-\infty, 0]} = \inf\{n\geq 1: S_n \leq 0\}
\]
 (where $(\alpha, \beta]$ denotes  the interval $\alpha <x\leq \beta$ as usual) and
define 
 \beqn\label{H/infty}
 H_{(-\infty,0]}^x(y)=P_x[ S_{T}=y],
 \eeqn
  the hitting distribution of $(-\infty, 0]$ for the walk starting at $x\in \mathbb{Z}$.  Likewise let $H_B^x$  be the hitting distribution of a non-empty set $B\subset \mathbb{R}$. [Thus $H_{(-\infty,0]}^1(y)=P[\hat Z=y-1], y\leq 0$ and  $H_{[0,\infty)}^{-1}(y)=P[ Z=y+1], y\geq 0$.]
   There exists $\lim_{x\to\infty} H_{(-\infty,0]}^x(y)$, which we denote by $H^{+\infty}_{(-\infty,0]}(y)$ and similarly for
   $H^{-\infty}_{[0,\infty)}$.  $H^{-\infty}_{[0,\infty)}$ is a probability distribution if  $E Z <\infty$ and vanishes identically otherwise \cite[P24.7]{S}.
  Let  $V_{{\rm ds}}(x), x =0, 1,  2,\ldots, $  be the renewal function  of the weak descending ladder-height process  (see (\ref{fr}) or  Appendix). For our present purpose it is convenient to bring in the function  $f_r$, the shift of $V_{{\rm ds}}$ to the right by 1, namely 
  \[
  f_r(x) = V_{{\rm ds}}(x-1)\qquad (x\geq 1).
  \]
    According to \cite[P19.5, E27.3]{S}, \cite[Section XII.3]{F}  $f_r$ is  a positive harmonic function on $[1,\infty)$, i.e.,  
    a positive solution of the equation 
    $f_r(x)= E_x[f_r(S_1);  S_1 \geq 1]$, which may be written as
 \beqn\label{H(z)}
    f_r(x) =\sum_{y= 1}^\infty f_r(y)p(y-x),\quad x
    \geq 1,
    \eeqn
and the solution  is unique apart from a constant factor;  it turns out  that the distribution of $Z$ is expressed as 
  \beqn\label{H(y)}
{P[Z > x]}= \sum_{y=1}^\infty f_r(y)p(y+x) \quad(x\geq 0),
\eeqn
  (see Theorem A  and (\ref{f_Z})  in Section 3 for more details).
 Define for any non-negative function $\fa (y)$, $y\leq 0$,
 \[
 H_{(-\infty,0]}^x\{\fa\}= E_x[\fa (S_T)]=\sum _{y\leq 0} H_{(-\infty,0]}^x(y) \fa (y) \leq \infty.
 \]
  
  For a set $B\subset \mathbb{R}$ such that $B\cap\mathbb{Z} \neq \emptyset$  let $g_B(x,y)$ denote the Green function of the walk killed as it hits $B$: 
\beqn\label{GF}
g_B(x,y)=E_x\Big[ \sum_{0\leq n<\sigma_B} \de(S_n,y)\Big] \quad\; \; (x, y\in \mathbb{Z}),
\eeqn
where $\de(x,y) =1$ if $x=y$  and $=0$  otherwise. This definition is 
different from that in \cite{S}, where the corresponding one agrees with our $g_B(x,y)$ if $x \notin B$,  but vanishes if $x\in B$ whereas according to our definition
\[
g_B(x,y) = \sum_{z\notin B}p(z-x)g_B(z,y) + \de(x,y)\quad \mbox{for}\quad x\in B, y\in \mathbb{Z}
\]
 (valid also for $x\notin B$); in particular $g_B(x,y)=\de(x,y)$ whenever $y \in B$. This relation  shows that  $g_B(x,y) $ equals the hitting distribution of $B$ by the dual (or time-reversed) walk started at $y$ which fact  is expressed as  
\[
g_B(x,y) = P_{-y}[ S_{\bar \sigma_{-B}}= -x]  \quad \mbox{for}\quad x\in B.
\] 
Here $-B=\{-z: z\in B\}$ and $\bar \sigma_B = \sigma_B$ if $S_0\notin B$ and $\bar \sigma_B =0$ otherwise.  

In case $B=(-\infty,0]$,  $g_B(x,y), x,y\in B$ is expressed explicitly  by means of the renewal functions of ascending and descending ladder height processes  (cf. Theorem A in Section 3),  by which it follows  immediately 
 that there exists $\lim_{y\to \infty} g_{(-\infty,0]}(x,y) $  which is denoted by  $g_{(-\infty,0]}(x,\infty) $ and  given by
\beqn\label{lim_g}
g_{(-\infty,0]}(x,\infty) = \left\{\begin{array}{ll} f_r(x)/EZ   \quad & x >0,\\
H^{-\infty}_{[0,\infty)}(-x) \quad & x \leq 0;
\end{array} \right.
\eeqn
if $EZ=\infty$,  the RHS vanishes so that $g_{(-\infty,0]}(x,\infty) =0$ for all $x$ (cf. (\ref{H-})).

\begin{theorem}\label{thm1.1} \, 
 {\rm (i)} \, For all  $x, y\in \mathbb{Z}$,
  \beqn\label{eqPS}
g_{(-\infty,0]}(x,y) + a(x-y) - H_{(-\infty,0]}^x\{a(\cdot-y)\} = A g_{(-\infty,0]}(x,\infty),
\eeqn
where   
\[
A=\left\{ \begin{array}{ll} 1/2 \quad &\mbox{if}\quad \sigma^2<\infty,\\
1\quad  &\mbox{if}\quad \sigma^2=\infty.
\end{array}\right.
\] 

{\rm (ii)}\;  If $EZ<\infty$, then as  $x\to\infty$, $a(x)/f_r(x) \to  A/EZ$ and  $a(-x)/a(x) \to 0$,  and 
 \[
 {\displaystyle  \sum_{x=0}^\infty a(-x)  P[\,  |X|> x] <\infty.}
 \]
 \indent
 {\rm (iii)}\; If $E Z=\infty$, then $\liminf_{x\to\infty} \,a(x)/f_r(x) = 0$.
\end{theorem}

   It is natural to extend $f_r(x)$ to a function on $\mathbb{Z}$ by means of (\ref{H(z)}) (so as to make  (\ref{H(z)}) valid for all $x\in \mathbb{Z}$), or what amounts to the same thing (in view of (\ref{H(y)})),
\beqn\label{f_ext}
f_r(x) = P[Z>-x]  \qquad \mbox{for}\quad  x\leq 0.
\eeqn
Since $f_r(0)=1< V_{\rm ds}(0)=f_r(1)$, $f_r$ is increasing. According to this \textit{extension} of $f_r$  together with the identity  
\beqn\label{eq/F}
H^{-\infty}_{[0,\infty)}(x) = P[Z>x]/E Z \quad (x\geq 0)
\eeqn
(cf. \cite[(XI.4.10)]{F} or the remark following (\ref{f_Z}))  relation  (\ref{lim_g}) is expressed simply  as 
 \beqn\label{g/fr}
 g_{(-\infty,0]}(x,\infty) =f_r(x)/EZ\quad\;\; (x\in \mathbb{Z}).
 \eeqn
For any integer $k$ and for any non-negative $\varphi$,  
$$H^x_{(-\infty,-k]}\{\varphi\} = H^{x+k}_{(-\infty,0]}\{\varphi(\cdot-k)\} \quad \mbox{and} \quad g_{(-\infty,0]}(x,z) 
= g_{(-\infty,-k]}(x-k,z-k)$$   
so that (\ref{eqPS})  is rephrased as 
    \beqn\label{H/a/A}
    H^{x}_{(-\infty,-k]}\{a(\cdot -y\} = a(x-y) + g_{(-\infty,-k]}(x, y) - Af_r(x+k)/EZ. 
    \eeqn

 The following corollary will be often useful in application of Theorem \ref{thm1.1}. For brevity of expression we write
   \[a^\dagger(x) = a(x) + \de(x,0)
   \]
  so that  $g_{(-\infty,0]}(x,y) + a(x-y) = a^\dagger(x-y)$ for $y\leq 0$.

\begin{corollary}\label{cor1} ~
 {\rm (i)}\; Suppose   $E Z<\infty$.  Then 
$H_{(-\infty,0]}^x\{a\}/ a(x)\to 0$ as $x\to\infty$ and 
\beqn\label{eq_Cor1}
 H_{(-\infty,0]}^x\{a(\cdot-y)\}= a^\dagger(x-y) -  A f_r(x)/ E Z \quad\mbox{for} \;\;\; x\in \mathbb{Z}, \,y\leq 0.
 \eeqn

 {\rm (ii)}\; If $E Z=\infty$, then 
  $H_{(-\infty,0]}^x\{a(\cdot-y)\}=a^\dagger (x-y)$ \;\;  for \;\;  $x\in \mathbb{Z}$, $y\leq 0$.
\end{corollary}

\begin{remark} \label{rem2.1}
(a)\, The statement  $a(-x)/a(x)\to 0$ in Theorem \ref{thm1.1} follows under the weaker condition $m_+(x)/m(x)\to 0$ ($x\to\infty$)  (cf. \cite[Theorem 4]{Upot}).  By Corollary \ref{cor1}(i) $a(x) \sim Af_r(x)/EZ$, which is shown in \cite{U1dm}  in case $\sigma^2<\infty$ and generalized   in \cite{Upot} 
 as $a(x)\sim f_r(x)/\int_0^xP[Z>t]dt$ under $m_+(x)/m(x)\to 0$  and $\sigma^2=\infty$.    
We include  proofs of these parts of Theorem \ref{thm1.1} which  are much simpler than  the proofs in \cite{Upot}---although the latter do not  depend on our Theorem \ref{thm1.1}.   
 
(b)\, By (\ref{H(y)}) it follows that 
\beqn \label{eq1.1}
\sum_{x=1}^\infty f_r(x)P[ X\geq x]=E Z.
\eeqn
This together with  Theorem \ref{thm1.1}(i) shows  that if $ E Z< \infty$,  then
 \beqn\label{a0}  
\sum_{x=1}^\infty [a(x)+a(-x)]P[ X> x] <\infty, \quad 
  \eeqn
 which  may be effectively used to derive  Chow's criterion for $EZ<\infty$  in a quite different way from   \cite{Cho} (see Remark \ref{rem3.1} for more details).
 
 (c)\,  
The process $M_n:=a(S_{n\wedge T})$ is a  non-negative martingale under $P_x$, $x\neq 0$, in particular $a(x)= E_x M_n$. Clearly $M_\infty = a(S_T)$ a.s.,  so that  $H_{(-\infty,0]}^x\{a\}=E_x M_\infty$. Hence
 Corollary  \ref{cor1} implies that   $(M_n)$ is uniformly integrable (so that $a(x)= E_x M_\infty$)   if and only if  $E Z=\infty$. 

 (d)\,  As another application of Theorem  \ref{thm1.1} we shall consider the random walk conditioned never to visit  the origin and observe that  the conditional walk  distinguishes $+\infty$ and $-\infty$ if and only if  either $E Z$ or $E \hat Z$ is finite, although its Martin compactification  does not  whenever  $\sigma^2=\infty$ (see Section 7).
 
  (e)\, Let $E Z=\infty$ and   consider asymptotic behaviour  of $a(x)/f_r(x)$ as $x\to\infty$.  Theorem \ref{thm1.1}(iii) tells  merely  $\liminf \, a(x)/f_r(x)= 0$. It however seems to be true quite generally that  $\lim a(x)/f_r(x) =0$. Actually if $E|X|<\infty$ and $\sigma_-^2:=E[X^2; X< 0]<\infty$  (in addition), one has $E|\hat Z|<\infty$ by virtue of the dual of 
(\ref{eq1.1}), so that  $f_r(x)\sim x/ E|\hat Z|$, hence 
  \beqn\label{a+a/f}
 [a(x)+a(-x)]/f_r(x) \longrightarrow 0 \qquad (x\to\infty).
 \eeqn
In case  $E|X|< \infty=\sigma_-^2$,  one has  $f_r(x) >\!> x/m_-(x)$ (cf. Lemma \ref{lem3.6}), and  if  (\ref{a/m})  is applicable,  this entails (\ref{a+a/f}). 
 If $X$ belongs to the domain of attraction of a stable law with exponent $1\leq \alpha \leq 2$ and skewness parameter $-1\leq \beta\leq 1$, then unless $\alpha=1$ and $\beta=0$, it follows that  $E|X|<\infty$  
and    $a(x)+a(-x) \sim \kappa x/m(x)$ for a positive constant $\kappa=\kappa_{\alpha, p}$ (cf. \cite[Section 8.1.1]{Upot}), hence (\ref{a+a/f}) holds; also in case $\alpha=1$ and $\beta=0$, if $ P[S_n>0]$ is further supposed to be convergent as $n\to \infty$,     (\ref{a+a/f}) holds (proved in Remark  \ref{rem5.1}).   
   \end{remark}

 The first part of the following theorem  provides asymptotic estimates of $a(x)$ as $|x|\to\infty$, and its third part   an answer to  the open question stated  at the very end of  Spitzer's book \cite{S} (see Remark \ref{rem2.2}(d) below).
 \begin{theorem}\label{thm1.2} \;  {\rm (i)}\, If $E Z<\infty$ and $ \sigma^2=\infty$,  
  then 
\beqn\label{c<c}  1  \leq  \liminf_{x\to\infty}  \frac{a(x)  m_-(x)}{x} \leq  \limsup_{x\to\infty}  \frac{a(x)  m_-(x)}{x} \leq 2;
\eeqn
and 
   \beqn\label{double}
\lim_{x\to\infty} \frac1{a(-x)} \sum_{z=1}^\infty P[z<Z\leq z+x]a(z) = EZ.
\eeqn
\vskip2mm\noindent
{\rm (ii)}\; If $E Z=\infty$, then $\lim_{x\to \infty} a(-x)=\infty$.
\vskip2mm\noindent
{\rm (iii)}\, Suppose  $\sigma^2=\infty =EZ$. Then  for some constant $C>1$
\beqn\label{a00}
C^{-1} a(-x) \leq \sum_{w=1}^x  \sum_{z=1}^\infty p(w+z) \Big[\frac{z}{m_-(z)}\Big]^2 \leq Ca(-x)  \quad (x\geq 1)  
\eeqn
[with  all the members vanishing if $P[X\geq2]=0$];
and there exists   $ \lim_{x\to\infty}a(-x) \leq \infty$  where the limit is finite if and only if  
\beqn\label{a01}
  \int_{1}^\infty \frac{t^2}{ m^2_-(t)}P[ X> t]dt<\infty
\eeqn
 and if this is the case,   $\lim_{x\to\infty}a(-x) = H^{-\infty}_{[0,\infty)}\{a\}$. 
\end{theorem}
 
\begin{remark}\label{rem2.2} (a)\,   (\ref{double}) entails that  $a(-x)$ is  asymptotically increasing   as $x\to\infty$  if  $ EZ <\infty$. Such  monotonicity of $a(-x)$  however  is verified under  $a(-x)/a(x) \to 0$  ($x\to\infty$) in  \cite[Corollary 40]{Upot} by a quite different approach.    
 
 (b)  In case $\sigma^2<\infty$  (\ref{a01})  holds whenever  $E[X_+^3]<\infty$ which condition is equivalent to $C^-:=\lim_{x\to -\infty} [a(x)+x/\sigma^2]<\infty$  \cite[Section 2.1]{U1dm}, while   (\ref{a01})  is possibly true even under the condition that for some 
$\delta>0$, $P[X>x] > x^{-1}(\log x)^{-1-\de}$ for all  sufficiently  large $x$. 

(c) \,Condition  (\ref{a01})  implies  $E Z<\infty$, the latter being equivalent to the integrability condition
  $\int_1^\infty t P[X> t]dt/m_-(t)<\infty$  (see \cite{Cho}, \cite[Section 2.4]{Unote}).
  
    (d)\,  
  Since  $\lim_{|x|\to\infty} a(x)=\infty$ if  $E|X|=\infty$ entailing $EX_+ = E X_- = \infty$ (because of the recurrence assumption), Theorem \ref{thm1.2}(iii) 
  gives an exact criterion for the trichotomy of   $\lim_{|x|\to\infty} a(x)=\infty$,  $M^-=\lim_{x\to -\infty} a(x)<\infty$,  $M^+= \lim_{x\to\infty} a(x)<\infty$. (This trichotomy itself is stated at the end of \cite{S}.)
  \end{remark}

  \begin{corollary}\label{cor2} Suppose  $\sigma^2=\infty$. There exists 
  $M^\pm:=\lim_{x\to \pm\infty}a(x) \leq \infty$ where  $M^- =0$  if and only if $P[X\geq 2]=0$ and
  in order that  $M^-<\infty$  each of the following conditions are necessary and sufficient.  
\vskip2mm\noindent  
 {\rm (i)} \quad  $\sum_{z=1}^\infty P[X>z]\big([a(z)]^2+a(-z)\big) <\infty$. 
\vskip2mm\noindent
 {\rm (ii)} \;  $ \sum_{z=1}^\infty P[Z>z]a(z)  <\infty$ and $P[X\leq -2]>0$.
\vskip2mm\noindent
  {\rm (iii)} \;  $ \sum_{z=1}^\infty  H_{[0,\infty)}^{x}(z)a(z)$ is bounded for  $x<0$ and $E|\hat Z| =\infty$.
\vskip2mm\noindent
 \end{corollary}
 \begin{proof}  The existence of the limit and the condition for $M^-=0$ follows immediately from
 Theorem \ref{thm1.2}.  Each of  conditions (i) and (ii)  implies $EZ<\infty$ (see Remark \ref{rem3.1}(b)  for (i) and note $M^+>0$  under  (ii)).  The assertion of the corollary then follows from Theorems \ref{thm1.1} and \ref{thm1.2} and the identity $H^{-\infty}_{[0,\infty)}(z) = P[Z>z]/E Z$. 
 \end{proof}

  For $y\in\mathbb{Z}$ write $\sigma_y$ for $\sigma_{\{y\}}$.   The  results (i) and (ii) given below are 
  taken    from Sections 7.3 and 7.5 of \cite{Upot}. 
  
{\rm  (i)} \;   If  $m_+(x)/m(x)\to 0$ ($x\to \infty$), then uniformly for $ 0\leq x\leq R$,
 as $R \to\infty$
\beqn\label{E11}
 P_x[\sigma_{[R,\infty)}< T] \sim P_x[\sigma_R < T]
 \sim \frac{f_r(x)}{f_r(R)},\;\mbox{and} 
\eeqn

{\rm (ii)} If   $m_+(x)/m(x)$  converges to 0 or  to  $1$ as $x\to\infty$, then  uniformly for $ 0\leq x\leq R$, as $R \to\infty$
\beqn\label{E0}
 P_x[\sigma_{[R,\infty)}<\sigma_0] \sim P_x[\sigma_R<\sigma_0]   \;\; \mbox{and}\;\; P_x[ T<\sigma_R] \sim P_x[\sigma_0 < \sigma_R]. \eeqn
 (The second relation of  (\ref{E0})  is the dual of  the first.)

These results are supplemented by   the following propositions. 

\begin{proposition}\label{prop1}
 \;
 If  $P[ X\geq 2] >0$, then for $x\in \mathbb{Z}$,
 \beqn\label{E2}
\lim_{R\to \infty} \frac{P_x[\sigma_R<T\,]}{P_x[\sigma_R<\sigma_0]} 
= \left\{\begin{array}{ll}
{Af_r(x)}/{[a^\dagger(x) EZ]} \quad  &\mbox{if} \quad EZ<\infty,\\
0 \quad  &\mbox{if} \quad EZ=\infty;
\end{array}\right.
\eeqn
where in case $EZ<\infty$  the convergence in (\ref{E2})   is uniform for $0\leq x\leq R$.
\end{proposition}
\begin{proposition}\label{prop2}\;  If $EZ<\infty$, then for $x\geq 0$,   as  $R-x\to \infty$,
\beqn\label{second}
P_x\big[\sigma_0<\sigma_{[R,\infty)}\big] \sim  P_x\big[ T <\sigma_{[R,\infty)}\big]  \sim \frac{f_r(R)- f_r(x)}{f_r(R)}\leq  \frac{f_r(R-x)}{f_r(R)}.
\eeqn 
\end{proposition}

The second equivalence in  (\ref{second}) follows from (\ref{E11})   if $x$ ranges  over  a set depending on $R$  in which  $f_r(x)=O(f_r(R)- f_r(x))$ but  does not otherwise.  If $P[ X\geq 2] =0$,   then ${P_x[\sigma_R<T\,]} ={P_x[\sigma_R<\sigma_0]}$ and $Af_r(x)/EZ=a^\dagger(x)$ for all  $x\in \mathbb{Z}$, where all the terms vanish for $x\leq -1$, and  the formula  (\ref{E2}) (necessarily the first case) is still
 reasonable. The dual statement  of  (\ref{E2}) for the case  $EZ<\infty$  may be written as
 \[
 \frac{P_x[\sigma_0<\sigma_{[R,\infty)}\,]}{P_x[\sigma_0<\sigma_R]} 
\sim  
\frac{Af_l(R-x)}{a^\dagger(-R+x) E|\hat Z|} \quad\mbox{uniformly for}\;\; 0\leq x\leq R\;\;  \mbox{if} \;\; 1< E|\hat Z|<\infty.
\]
The corresponding one for (\ref{second}) will be stated as Lemma \ref{lem6.3} in Section 6.

  The formula  (\ref{E2})
   says that if $EZ<\infty$, then  $ P_x[\sigma_R< T\,|\, \sigma_R<\sigma_0]$, being  equal to  the ratio on the LHS, approaches  unity as $x$ becomes  large independently of how $R$ is large, 
 while   if $EZ=\infty$ this is not the case:  this conditional probability tends to zero as $R\to\infty$,
 in other words,  for $R$ large enough the walk---even if it is conditioned  
 on $\sigma_R<\sigma_0$---reaches $R$    only  after entering the negative half line with
overwhelming  probability as far as  its starting position  $x$ is fixed.
If $m_+/m\to 0$ and $\ell_+(x)= \int_0^x P[Z>t]dt$, then $P_x[\sigma_R< \sigma_0\,] \sim a^\dagger(x)/a(R)$ (see Lemma \ref{lem6.1}),  $\ell_+$ is slowly varying and $f_r(x)\sim a(x)\ell_+(x)$ \cite[Lemma 46]{Upot}, so that  by (\ref{E11}) it follows that as $x\to \infty$ under
 $x<R$
\beqn\label{Af/a} 
 P_x\big[\sigma_R <T \,\big|\, \sigma_R<\sigma_0\big] \sim \ell_+(x)/\ell_+(R).
  \eeqn
This shows that for each $\ep>0$, the ratio above approaches 1 as $R\to\infty$ uniformly for $x>\ep R$.  The same holds true if $E|\hat Z|<\infty$ at least under some regularity condition on the tails of $F$  but can fail  in general (see Remark \ref{rem6.2} of Section 6).

The rest of the paper is organized as follows.
In Section 3 we collect  fundamental facts used in this paper about $f_r$, $a(x)$, $g_{(-\infty,0]}$ etc.
given in Spitzer \cite{S} and advance several  lemmas that are directly derived from them.
The proofs of Theorems \ref{thm1.1} and \ref{thm1.2} are given in  Sections 4 and 5, respectively. The proofs of  Propositions \ref{prop1} and \ref{prop2}  are given in  Section 6.
 In Section 7 we briefly study    large time behaviour of the walk conditioned never  to visit the origin.   In Section 8 (Appendix)   we present   a few facts about strictly  and weakly ascending ladder height variables.


\section{ Preliminary lemmas }

In this section we collect fundamental results of the recurrent random walks on $\mathbb{Z}$ given in Spitzer's book \cite{S} and  then derive  some  
consequences of them that are used the later sections. 

 For $B \subset \mathbb{Z}$ we have defined the first hitting time by $\sigma_B=\inf\{n\geq 1: S_n\in B\}$.  For a point  $x\in \mathbb{Z}$ write $\sigma_x$ for $\sigma_{\{x\}}$.  For typographical reason we sometimes write
$\sigma B$ for $\sigma_B$.   

Let  $u_{\rm as}(x)$, $x=0, 1, 2,\ldots$ be the renewal sequence of the strictly ascending 
ladder variables, namely  $u_{\rm as}(0)=1$ and
\beqn\label{v+}
u_{\rm as}(x)= \sum_{n=1}^\infty P[ Z_1+\cdots + Z_n = x]\qquad x\geq 1;
\eeqn
and similarly  $v_{\rm ds}(x)$, $x=0, 1, 2,\ldots$   denotes  the renewal sequence of the \textit{weak descending}  
ladder variables, which may be given   by $v_{\rm ds}(0) =1/c$ and
\beqn\label{v-}
v_{\rm ds}(x)= \frac1{c} \sum_{n=1}^\infty P[ \hat Z_1+\cdots  + \hat Z_n = - x]\qquad x\geq 1,
\eeqn
where
\[
c=\exp\left[ -\sum_{k=1}^\infty \frac1{k} p^k(0)\right]=\exp\left[\frac1{2\pi}\int_{-\pi}^\pi \log \big|1-E[e^{itX}] \big| dt \right]\!.
\]
  (See Appendix  for  (\ref{v-}) as well as for  the probabilistic meaning of  the constant $c$.)
Owing to the renewal theorem \cite{F},  there exist  limits
\beqn\label{v(+infty)}
u_{\rm as}(\infty):=\lim_{x\to\infty}u_{\rm as}(x)=1/EZ \qquad\mbox{and}\qquad v_{\rm ds}(\infty):=\lim_{x\to\infty} v_{\rm ds}(x)=1/ c\, E[- \hat Z]\,.
\eeqn
 The Green function   $g_{B}(x,y)$ ($\,x,y \in \mathbb{Z}\,$)  defined in (\ref{GF}) may be written as: 
\[
g_{B}(x,y)=\sum_{n=0}^\infty P_x[S_n=y, n< \sigma_B ].
\]
The following theorem  follows from the propositions P18.8, P19.3, P19.5 of \cite{S}.  
For two real numbers $s$ and $t$  write $s\wedge t =\min\{s,t\}$ and  $s\vee t=\max \{s,t\}$.

\vskip2mm\noindent
{\bf Theorem A.} \; \;  
(i) \quad $u_{\rm as}(\infty)v_{\rm ds}(\infty) =1/{c EZE|\hat Z|} =2/\sigma^2$.
\vskip3mm

(ii) \quad  $g_{\,(-\infty,0]}(x,y)= \sum_{z=1}^{x\wedge y} v_{\rm ds}(x-z)u_{\rm as}(y-z)\qquad (x,y >0) $; \;\;\textit{and }
\vskip3mm
  
 \qquad\;
${\textstyle  g_{\,[0,+\infty)}(x,y)= g_{\,(-\infty,0]}(-y,-x)= \sum_{z=1}^{|x|\wedge|y|} u_{\rm as}(|x|-z)v_{\rm ds}(|y|-z)\qquad (x,y <0).}$

\vskip2mm\vskip2mm

The formulae in Theorem A will  often be  used in combination  with the following representation of the hitting distribution 
 $H_{(-\infty,0]}^x(y)$ of $(-\infty,0]$:
 \beqn\label{H}
 H_{(-\infty,0]}^x(y)= \sum_{z=1}^\infty g_{(-\infty,0]}(x,z)p(y-z)\qquad (x>0, y\leq 0),
 \eeqn
 and analogous one for $H_{[0,\infty)}^x$ (see (\ref{rep_H}) for another representation). The function $f_r$ may be written as 
\beqn\label{fr}
f_r(x)=v_{\rm ds}(0)+\cdots +v_{\rm ds}(x-1)  \qquad\quad (x\geq1),
\eeqn
and its  dual   as\,
  $f_l(x)=c^{-1}\,[u_{\rm as}(0)+\cdots +u_{\rm as}(x-1)] \; (x\geq1)$.

By Theorem A(ii)  and $u_{\rm as}(y)\leq 1$ it follows that  
 \beqn\label{eqTB}
 g_{(-\infty,0]} (x,y) \leq  \left\{   \begin{array}{ll}  f_r(x) \quad &\mbox{if} \;\, x\leq y,\\
 f_r(x)-f_r(x-y) \quad &\mbox{if}\;\; x> y.
\end{array}\right.
 \eeqn
  Let $x\to -\infty$ in   
  $H_{[0,\infty)}^x(y) = \sum_{z<0} g_{[0,+\infty)} (x,z) p(y-z)$.  Noting  $g_{[0,+\infty)} (x,z) = g_{(\infty,0]} (-z, -x)   \to f_r(-z)/EZ$ and $\sum_{w\geq1} f_r(w)p(y+w) <\infty$,  we then find  that 
 \beqn\label{H-}
H_{[0,\infty)}^{-\infty}(y):=\lim_{x\to\, -\infty}H_{[0,\infty)}^x(y) =\frac1{E Z} \sum_{w=1}^\infty f_r(w) p(y+w)
  \eeqn
It also follows that  $H_{[0,\infty)}^x(y)  \leq \sum_{w=1}^\infty f_r(w) p(y+w) $, so that
\beqn\label{H<H}
H_{[0,\infty)}^x(y)  \leq  (E Z) H_{[0,\infty)}^{-\infty}(y) \quad\mbox{for all} \quad x\leq 0<y \quad\mbox{if $EZ <\infty$.}
\eeqn
In particular  the  three conditions
 ${\rm (a)}  \;  E Z = \infty; \; {\rm (b)}\;  u_{\rm as}(\infty)=0; \; {\rm (c)} \;  H_{[0,\infty)}^{-\infty}(\cdot) \equiv 0$
  are equivalent to one another.
 Since  $g_{[1,\infty)}(0,-y)=g_{(-\infty,0]}(y+1, 1) =v_{\rm ds}(y)$  we have for  $k>0$
\beqn\label{sharp}
P[Z=k] =\sum_{y=0}^\infty g_{[1,\infty)}(0,-y)p(k+y) =  \sum_{y=0}^\infty v_{\rm ds}(y)p(k+y), 
\eeqn    
and, by  summation by parts, 
\beqn\label{f_Z}
P[Z  >  x] = \sum_{y=0}^\infty v_{\rm ds}(y)P[X > x+y] 
= \sum_{y=1}^\infty f_r(y)p(x+y) \qquad (x\geq 0).
\eeqn
Note that  the last equality together with (\ref{H-})  yields (\ref{eq/F})  (i.e., $H_{[0,\infty)}^{-\infty}(x)=P[Z>x]/EZ$).

The next theorem also is taken from Spitzer \cite[T28.1, T29.1, P30.2, P30.3]{S}.
\vskip2mm\noindent
{\bf Theorem B.}  \textit{ The series  $\sum_{n=0}^\infty[p^n(0)-p^n(-x)]$ converges for each $x\in \mathbb{Z}$ and if $a(x)$ denotes the sum, then  
 the following relations hold.
\beqn\label{g}
g_{\{0\}}(x,y)=  a^\dagger (x)+a(-y)-a(x-y) \qquad (x, y \in \mathbb{Z}),
\eeqn
 \beqn\label{s_add/a}
 a(x+y)\leq a(x)+a(y) \quad\mbox{and}\quad   a^\dagger(x)+ a(-x) \geq 1 \qquad (x, y \in \mathbb{Z}),
 \eeqn
\beqn \label{eq2.6}
\sum_{z=-\infty}^\infty p(z-x)a(z-y) =a^\dagger(x-y),
\eeqn  
\beqn\label{a1}
\lim_{x\to \pm \infty} [a(x+1)-a(x)] = \pm 1/\sigma^2\quad \mbox{and}\quad \lim_{x\to\infty} [ a(x)+a(-x)]=\infty.
\eeqn  
 If the walk is left-continuous (i.e. $P[ X\leq -2]=0$), then $a(x)=x/\sigma^2$ for $x>0$; analogously  $a(x)=-x/\sigma^2$ for $x<0$ for right-continuous walks;  except for  left- or right-continuous walks with infinite variance $a(x)>0$ for all $x\neq 0$.}

[(\ref{g}) with $x=0$ and the second inequality of (\ref{s_add/a}), not given in \cite{S}, follows from (\ref{eq2.6}) and $g_{\{0\}}(x,x)\geq 1$, respectively.] 

\vskip2mm
 We put
\[
\bar a(x)=\frac12 [a(x)+a(-x)].
\]
By (\ref{g})  it follows  that 
$g_{\{0\}}(y,y) = 2 \bar a(y) +\de(0,y)>0$
  and that 
\beqn\label{a3}
P_x[\sigma_y<\sigma_0]= \frac{g_{\{0\}}(x,y)}{g_{\{0\}}(y,y)} =\frac{a^\dagger(x)+a(-y)-a(x-y)}{2\bar a(y)} \qquad  (x, y\in \mathbb{Z}, y\notin \{x,0\}).
\eeqn
The  equation (\ref{eq2.6}) states that $a(x)$ is harmonic on $x\neq 0$, which together with $a(0)=0$  entails that the process $M_n:=a(S_{\sigma_\xi \wedge n} -\xi)$ is a martingale, provided that $S_0\neq \xi \in \mathbb{Z}$ a.s.  Using  the optional sampling theorem and  Fatou's lemma  we obtain first  the inequality
 $a(x-\xi) = \lim_{n\to\infty} E_x[M_n] \geq E_x[a(S_{\sigma_\xi \wedge\sigma_B}-\xi)] $
valid whenever  $ x\neq \xi$, and then by using (\ref{eq2.6}) again if $\xi\in B$, 
\beqn\label{aa0}
E_\xi[a(S_{\sigma_B}-\xi)] =\sum_{y\in B} p(y-\xi) a(y-\xi)+  \sum_{z\notin B}p(z-\xi) E_z[a(S_{\sigma_\xi \wedge\sigma_B}-\xi)] \leq 1,
\eeqn
 so that 
\beqn\label{aa1}
E_x[a(S_{\sigma_B}-\xi)] \leq a^\dagger(x-\xi) \quad \mbox{for}\quad \xi\in B, x\in \mathbb{Z}, 
\eeqn
 in particular  
\beqn\label{eq2.-1}
 a(y)P_x[\sigma_y<\sigma_0] = E_x[a(S_{\sigma_0\wedge \sigma_y})] \leq a^\dagger (x) \qquad  (x, y\in \mathbb{Z}).
\eeqn

  In the rest of this section we prove several lemmas that  are derived more or less directly from the results presented above.   

\begin{lemma}\label{lem3.1}\, Let  $\sigma^2=\infty$.  Then there exists $\lim_{x\to \infty}a(x) ~(\leq \infty)$ which  is zero if
and only if the random walk is left-continuous. 
\end{lemma}
By the last statement of Theorem B this lemma  shows that 
$\inf_{x\geq 1} a(x) >0$ unless  the random walk is left continuous. 

\begin{proof}  Let $\sigma^2=\infty$.  The relations  (\ref{eq2.-1}) and (\ref{a3})  yield
\beqn\label{L3.1}
a(y)\geq \frac{a(x)}{a(x)+a(-x)}\big[a(y)+a(-x)-a(y-x)\big] \quad \quad (x\neq 0).
\eeqn
 On using  (\ref{a1})  it then follows that
\beqn\label{eq_lm_2.4.2}
\liminf_{y\to\infty} a(y) \geq \,  \frac{a(x)a(-x)}{a(x)+a(-x)}   \quad \mbox{for all}\quad x\neq 0.
\eeqn
If $\limsup_{x\to\infty} a(x) <\infty$, then $\lim_{x\to\infty} a(-x) =\infty$  in view of (\ref{a1}) and the inequality (\ref{eq_lm_2.4.2}) gives $\liminf a(x)\geq \limsup a(x)$ so that $\lim a(x)$ exists.  If this limit is zero, then the RHS of (\ref{eq_lm_2.4.2}) must be zero for all $x>0$, which  is possible only if the walk is left-continuous.

Now suppose $\limsup_{x\to\infty} a(x) = \infty$ and  put  $M=\liminf_{x\to\infty} a(x) (\leq \infty)$.  Contrary to what is to be shown let $M<\infty$.  Then one can  choose $R$ such that $a(x)+a(-x)>4M+6$ for $x >  R$. In view of (\ref{a1}) there must exist $x_1>R$ such that $2M+2\leq a(x_1)< 2M+3$, which entails $a(-x_1)> 2M+3$. Combined with (\ref{eq_lm_2.4.2}) these lead to the absurdity
\[
M\geq  \frac{a(x_1)a(-x_1)}{a(x_1)+a(-x_1)} \geq\frac{a(x_1)}{2}\geq M+1.
\]
Hence  $M$ must be infinite. 
\end{proof}

\begin{lemma}\label{lem3.2}  For all  $x, y\in \mathbb{Z}$,
\beqn\label{L6.1}
-\, \frac{a(y)}{a(-y)} a(x)\leq a(x+y)-a(y) \leq a(x) \quad \mbox{if} \;\; a(-y) \neq 0.
\eeqn
\end{lemma}

\begin{proof}    From  (\ref{a3}) and (\ref{eq2.-1}) we have
\[
\frac{a(x) + a(y) - a(x+y)}{a(y)+a(-y)} \leq \frac{a(x)}{a(-y)} \quad (a(-y)\neq 0)
\]
 [in (\ref{L3.1}) take $-y$ and $x$ in place of $x$ and $y$ respectively],  which, after simple rearrangements,  becomes the left-hand inequality of (\ref{L6.1}); the case $x=0$ is obvious.
The right-hand one is the same as $g_{\{0\}}(x,-y)\geq 0$.  \end{proof} 
  
 Put $g(x,y)=g_{\{0\}}(x,y) -\de(x,0)$, or explicitly $g(x,y) = a(x)+a(-y)-a(x-y)$. 
\begin{lemma}\label{lem3.3} 
If $B$ is a  proper subset of\,  $\mathbb{Z}$ such that  $0\in B$, then 
\beqn\label{g/g/a}
g_{\{0\}}(x,y)  = g_B(x,y) +E_x[g(S_{\sigma_B}, y)]\quad (x, y \in \mathbb{Z}).
\eeqn 
\begin{proof}  Let $\La_B(y)$ be the number of  visits  to $y$ in the time interval
 $\{1,2,\ldots,\sigma_B-1\}$: 
\[ \La_B(y) = \sharp\{n\geq 1: S_n=y, n<\sigma_B\}. 
\]
Then $g_{\{0\}}(x,y) =\de(x,y) + E_x[\La_{\{0\}}(y)]$ and similarly for $g_B(x,y)$, and  (\ref{g/g/a}) can be written as
\beqn\label{La/g}
E_x[\La_{\{0\}}(y)] = E_x[\La_{B}(y)]  + E_x[g(S_{\sigma_B}, y)],
\eeqn
 provided that $0\in B$ which
entails $\sigma_B\leq \sigma_0$ a.s. 
 Recall that $g(0,y)=0$ and for $z\neq 0$, $g(z,y) =g_{\{0\}}(x,y)$, the expected number 
of  visits to $y$ before entering.   
 If $y\notin B$, then by the strong Markov property the
above equality follows immediately,  It therefore suffices to show (\ref{La/g}) for $y\in B$.

Let $y\in B$, when one always has $\La_B(y) =0$ a.s. For $x\notin B$,  (\ref{La/g})  then follows   immediately. For $x\in B$, one observes that
\[E_x[g(S_{\sigma_B}, y)] =\sum_{z\in B} p(z-x)g(z,y) + \sum_{z\notin B} p(z-x)E_z[g(S_{\sigma_B},y)],
 \]
but (\ref{La/g}) with  $x$ replaced by  $z\notin B$ is valid  so that $E_z[g(S_{\sigma_B},y)] = E_z[\La_{\{0\}}(y)]$, of which the RHS equals   $g(z,y)$ for $z\neq 0$.   Thus 
\[E_x[g(S_{\sigma_B}, y)] =\sum_{z\neq  0}p(z-x)g(z,y),
\]
which shows  (\ref{La/g}), for  the last sum equals $E_x[\La_{\{0\}}(y)]$.
This finishes the proof.  
\end{proof}

\end{lemma} 

\begin{lemma}\label{lem3.4} If the walk is not left-continuous and  $k^+:=\sup_{x\geq 1} a(-x)/a(x)$, then  
\beqn\label{eL3.3}
0\leq g_{\{0\}}(x,y) - g_{(-\infty,0]}(x,y) \leq (1+k^+)a(-y) \quad(x,  y \in \mathbb{Z}). 
\eeqn
  If the walk is not right-continuous and  $k^-:=\sup_{x\leq -1} a(-x)/a(x)$, then  
\beqn\label{eL3.31}
0\leq g_{\{0\}}(x,y) - g_{(-\infty,0]}(x,y) \leq (1+k^-)a^\dagger(x) \quad(x. y \in \mathbb{Z}). 
\eeqn
\end{lemma}
\begin{proof}    
Take $B=(-\infty, 0]$ in  (\ref{g/g/a}) and use 
the inequality $a(z)-a(z-y)\leq \big[a(z)/a(-z)\big]a(-y)$ ($z\leq-1$) that follows from (\ref{L6.1}) to 
see that the difference on the middle member of (\ref{eL3.3}) is  not larger than
\[
E[g(S^{x}_{\sigma(-\infty,0]}, y)]  \leq (1+k^+)a(-y),
\]
hence the right-hand inequality of (\ref{eL3.3}). The left-hand one is trivial.

The  right-hand inequality of (\ref{eL3.31}) is also derived from (\ref{g/g/a}) but this time we  use the inequality $g(z, y) \leq 2\bar a(z)$  to have
\beqn\label{1+k}
E_x[g(S_T, y)  ]  \leq 2E_x[\bar a(S_{T})] =2H^{x}_{(-\infty,0]}\{\bar a\}.
\eeqn
By definition 
$2\bar a(z) \leq (1+k^-)a(z)$  for $z\leq 0$, while  $H_{(-\infty,0]}^{x}\{a\}\leq a^\dagger(x)$ as a special case of  (\ref{aa1}). Hence 
$2H_{(-\infty,0]}^{x}\{\bar a\} \leq  (1+k^-)a^\dagger(x)$, showing  (\ref{eL3.31}). 
 \end{proof}
 
In \cite{Upot},  Lemma \ref{lem3.4} plays a significant role for the proof of (\ref{E11}). In this article we apply it only to obtain the next result. 
 
\begin{lemma}\label{lem3.5}\,  If either $a(-x)/a(x)\to 0$  or $a(x)/a(-x)\to 0$   as $x\to\infty$ (with the understanding that $a(x) >0$ ($a(-x)>0$) for $x>0$ in the former (latter) case),  then
\beqn\label{2.7}
\lim_{x\to\infty}\frac{g_{(-\infty,0]}(x,x)}{g_{\{0\}}(x,x)}= 1;
\eeqn
and
\beqn\label{f/a}
 \left\{\begin{array}{ll}
 \lim_{x\to\infty}\, f_r(x)/2\bar a(x)= E Z \leq \infty \quad &\mbox{ if\;\;  $ a(-x)/a(x) \to 0$,}\\
\lim_{x\to\infty}\, {f_l(x)}/{2\bar a(x)}= -E \hat Z \leq \infty \quad &\mbox{ if\;\;  $a(x)/a(-x) \to 0$.}
\end{array}\right.
\eeqn
The identities in (\ref{2.7}) and (\ref{f/a}) are valid whenever $\sigma^2<\infty$.
\end{lemma}
\begin{proof}  Under the assumption of the lemma  it follows from Lemma \ref{lem3.4} that  $g_{[0,\infty)}(x,x) = g_{\{0\}}(x,x) + o(\bar a(x)) = 2\bar a(x)\{1+o(1)\}$. This verifies  (\ref{2.7}),   which  entails (\ref{f/a})   in view of  Theorem A(ii). In case  $\sigma^2<\infty$  use  the explicit asymptotic forms of $v_{\rm ds}$, $u_{\rm as}$ and $\bar a$ to deduce  (\ref{2.7}) from Theorems A and B; then observe  $g_{-\infty,0]}(x,x) \sim f_r(x)/EZ$ to see  $f_r(x)/EZ \sim 2\bar a(x)$.
\end{proof}

\begin{lemma}\label{lem3.6} \;\;  If $ \sigma_-^2 := E[  X_-^2] =\infty> E X_-$,  then 
\vskip2mm\noindent
{\rm (i)}  \qquad ${\displaystyle 
\frac1{m_-(x)}\int_0^x P[\hat Z<-t]dt ~\longrightarrow ~\frac1{c E Z} 
}$~~ as~ $x\to\infty$;\,\mbox{and}
\vskip2mm\noindent
{\rm (ii)} \qquad ${\displaystyle 
E Z\leq  \liminf_{x\to\infty}  \frac{f_r(x)  m_-(x)}{x} \leq  \limsup_{x\to\infty}  \frac{f_r(x)  m_-(x)}{x} \leq 2E Z.
}$
\vskip2mm\noindent
 If $ \sigma_-^2 <\infty$,  the   $\liminf$ and $\limsup$   in {\rm (ii)} coincide and equal $m_-(+\infty)/cE|\hat Z| \in [0,\infty)$.
\end{lemma}
\begin{proof}     As a dual relation of  (\ref{f_Z})  we have for $t\geq 0$
\beqn \label{hat Z}
P[\hat Z<-t] = v_{\rm ds}(0) \sum_{y=0}^\infty u_{\rm as}(y)P[ X<-t-y].
\eeqn
Let $\sigma_-^2=\infty> EX_-$, which  entails that 
$\int_0^x P[  X< -t - y] dt/ m_-(x)$ tends to zero as $x\to\infty$ for each $y\geq 0$.   Replacing
 $u_{\rm as}(y)$ by  $u_{\rm as}(\infty) +o(1)$ in (\ref{hat Z})  and recalling  $v_{\rm ds}(0)u_{\rm as}(\infty) =1/c E Z$  we then  infer that
\[
\frac1{m_-(x)} \sum_{k=0}^x P[\hat Z<-k] =\frac{v_{\rm ds}(0)}{m_-(x)} \sum_{k=0}^x \sum_{y= k}^\infty u_{\rm as}(y-k)P[ X<-y]  =  \frac1{c E Z}+ o(1).
\]
Thus (i) is verified. Noting that $c f_r(x+1)$ is the renewal function for the variable $-\hat Z$ we use
 the first inequality of Lemma  1 of Erickson  \cite{E} which may read
\[
1\leq \frac{c f_r(x+1)}{x} \int_0^x P[\hat Z<-t]dt \leq 2;
\]
combining this with (i)  we can readily deduce (ii). The last assertion is obvious, for  $ m_-(\infty)<\infty$ if $\sigma^2_-<\infty$     and $f_r(x)/x \to 1/c E|\hat Z|$.
\end{proof}

\begin{remark}\label{rem3.1}\, By (\ref{f_Z}) and Lemma \ref{lem3.6}(ii)  one infers that  $\int_0^\infty  P[X>t]tdt/m_-(t)<\infty$ if $EZ<\infty$---the necessity part of the Chow's criterion for $EZ <\infty$ (this half of it is also proved by Doney  \cite{D2}).  Combined with (\ref{a0}) this shows that if $EZ<\infty$, then  both of the following summability conditions
 hold 
\[  (\sharp)\; \; \sum \bar a(x) P[X>x] <\infty\;\;\;\;\;  \mbox{and} \;\;\;\;\;  (\flat)\;\;\int_1^\infty \frac{x P[X>x]}{m_-(x)}dx<\infty.
\]
The converse as well as the implication $(\flat) \Rightarrow (\sharp)$    is proved in \cite[Lemma 4.1, Lemma 2.9, Proposition 2.1(i)]{Unote}. Thus the equivalence of $(\flat)$ and $EZ<\infty$ follows.  It also holds that  $(\sharp)$ and $(\flat)$ are equivalent  \cite[Corollary 4.1]{Unote}, \cite[Eq(1.3)]{Upot} as  (\ref{a/m}) may suggest.
\end{remark} 
  
\begin{lemma}\label{lem3.7} \, Suppose $E Z<\infty$ and $\sigma^2=\infty$. Then  
$\lim_{x\to\infty}\, {a(-x)}/{ a(x)}= 0.$
\end{lemma}
\begin{proof}  
This lemma follows from Theorem 4 of \cite{Upot} as mentioned previously.  The direct proof is easy and given as follows. Let $\bar \sigma_x =\sigma_x$ if   $S_0\neq x$ and $=0$ otherwise. Then, 
applying (\ref{a1}) and  (\ref{H-})   under the assumption of the lemma leads to
\[
\frac{a(-x)}{a(x)+a(-x)}=\lim_{z\to \,-\infty}P_z[\sigma_x<\sigma_0] = \sum_{y=1}^\infty H_{[0,\infty)}^{-\infty}(y)P_y[\bar \sigma_x <\sigma_0].
\]
As $x\to\infty$ the last sum approaches zero and hence $a(-x)/a(x)\to 0$. 
\end{proof}

\begin{lemma}\label{lem3.8} \;  If $E Z<\infty$, then  
\[
 \sum_{y=-\infty}^0 a(y) P[X<y] <\infty  \quad\mbox{and}  \quad \lim_{x\to +\infty} \frac{H_{(-\infty,0]}^x\{a\}}{f_r(x)} =0.
\]
\end{lemma}
\begin{proof}  
By Theorem A(ii)  $g_{(-\infty,0]} (1,z) = v_{\rm ds}(0)u_{\rm as}(z-1)$ for $z\geq 1$.  Suppose  $E Z<\infty$. Then $u_{\rm as}(\infty) >0$  and for $y\leq 0$,
\[H_{(-\infty,0]}^1(y)=  v_{\rm ds}(0)\sum_{z=1}^\infty  u_{\rm as}(z-1)p(y-z) \asymp P[X<y],\]
and hence   the first assertion  follows,  for   $ H_{(-\infty,0]}^1\{a\} <\infty$ by virtue of  (\ref{aa1}).  Since $H_{(-\infty,0]}^x(y)$ is less than $f_r(x) P[ X<y]$ and $H_{(-\infty,0]}^x(y)/f_r(x)\to 0$ as $x\to\infty$ for each $y\leq 0$, by dominated convergence
$H_{(-\infty,0]}^x\{a\}/f_r(x)\to 0$, as desired. 
\end{proof}

In view of the following lemma  we can define for any non-empty  subset  $B$ of $ \mathbb{Z}$ the function  $u_B(x)$, $x\in \mathbb{Z}$ by 
 \beqn\label{eqPS1}
u_B(x)= g_{B}(x,y) +  a(x-y) - H_{B}^x\{a(\cdot-y)\}.
\eeqn

\begin{lemma}\label{lem3.9} \, For each   $x  \in \mathbb{Z}$ the RHS of (\ref{eqPS1})
is  independent of $y\in \mathbb{Z}$, and $u_B$ defined therein  is  non-negative,  represented by $u_B(x) = a^\dagger(x-\xi)-H^x\{a(\cdot -\xi)\}$ for any $\xi\in B$ and harmonic on $\mathbb{Z}\setminus B$  in the  sense  that
for each $\xi\in B$ fixed,  
 \beqn\label{eqPS5}
\sum_{z\notin B} p(z-x)u_B(z) = u_B(x) \quad\mbox{for}\quad x\in \mathbb{Z}.
\eeqn
\end{lemma}

Identity (\ref{eqPS1}) and hence what are advanced below hold true for every recurrent random walk irreducible on $\mathbb{Z}$.
The analogous result holds for the two-dimensional recurrent random walks to which the same proof  applies. 

\begin{proof}  
In the proof of Lemma 2.9 of \cite{Uf.s} it is shown that  for each $x\in \mathbb{Z}$ fixed, the RHS of   (\ref{eqPS1}) is  a dual-harmonic function of $y\in \mathbb{Z}$ (i.e., harmonic with respect to the dual transition function $\hat P(x,y): = p(x-y)$).   In view of the uniqueness theorem of  non-negative harmonic function \cite[P13.1]{S}, \cite[Proposition 6-3]{KSK} the first assertion of the lemma accordingly    follows if we show that   it is bounded below. 
To this end it suffices to see that for all $x, y\in \mathbb{Z}$ and $\xi\in B$,
\[
 H_{B}^x\{a(\cdot-y)\} \leq a^\dagger(x-\xi) + a(\xi-y)
\]
(since $a(x-y)- a(\xi -y)$ is a bounded function of $y$),  which however is  immediate from the  subadditivity
$a(\cdot-y) \leq a(\cdot-\xi)+ a( \xi-y)$ and the inequality (\ref{aa1}).
  Taking  $y$ from $B$  in (\ref{eqPS1}) it follows that
  \beqn\label{u/a/Ha} u_B(x) =a^\dagger(x-\xi) - H_B^x\{a(\cdot -\xi)\} \qquad \mbox{for}\;\; \xi\in B
  \eeqn
 and  by the inequality (\ref{aa1})  we  see    $u_B\geq 0$. 
  Noting that  
$ \sum_{w\notin B} p(w-x) H_{B}^w\{a(\cdot-\xi)\}=
H_{B}^x\{a(\cdot-\xi)\}-   \sum_{z\in B} p(z-x) \{a(z-\xi)\}$
 and using (\ref{eq2.6})  one deduces
\beqn\label{harm}
\sum_{w\notin B}p(w-x)\big[a(w-\xi)- H_{B}^w\{a(\cdot-\xi)\}\big] =a^\dagger (x-\xi)- H_{B}^x\{a(\cdot-\xi)\},
\eeqn
which shows  (\ref{eqPS5}), for $a(w-\xi)$ can be replaced by  $a^\dagger(w-\xi)$ because of the identity   $\sum_{w\notin B}p(w-x)\de(w,\xi)=0$.
 \end{proof}

\begin{remark}\label{rem3.2}  (a)\,  The independence of  the RHS of (\ref{eqPS1}) from $y$ also  follows from  Lemma \ref{lem3.3}. 
Indeed,  
 if $0\in B$ then we have  (\ref{g/g/a}) which  becomes (\ref{eqPS1}) with $u_B(x) = a^\dagger(x) -  H^x_B\{a\}$ after a simple rearrangement of terms.  For the case  $0\notin B$,   pick any $\xi\in B$,  consider  (\ref{g/g/a}) for $B'= B -\xi$ (shift by $\xi$)  in place of $B$ and replace $x$, $y$ by $x-\xi$ and $y-\xi$.  Conversely   Lemma \ref{lem3.3} follows immediately from the first half of Lemma \ref{lem3.9}.

(b) \, For a positive integer $R$ let $\tau_R = \sigma_{\mathbb{Z}\setminus (-R, R)}$. Then
\beqn\label{a/Ea}
a^\dagger(x-\xi)=   E_x[a(S_{\tau_R\wedge \sigma_B}-\xi)]   \qquad (x\in \mathbb{Z}, \xi\in B),
\eeqn
and  the function  $u_B$ defined in (\ref{eqPS1}) is given by 
\beqn\label{eq_auB}
u_B(x)=  \lim_{R\to\infty} E_x\big[a(S_{\tau_R}-\xi); \tau_R < \sigma_B\big]
\qquad (x\in \mathbb{Z},\xi\in B);
\eeqn
in particular the limit appearing in (\ref{eq_auB}) is independent of the choice of $\xi$.   
 These formulae are verified as follows.
 For $x\neq \xi$, $M_n := a(S_{n\wedge\sigma_B}-\xi)$ being a non-negative martingale under $P_x$ that is uniformly bounded on $n<\tau_R$, one obtains the identity  $a(x-\xi)=  E_x M_{\tau_R}$. As for the case $x=\xi$ suppose that   $\xi=0\in B$ for simplicity so that $M_{\tau_R}= a(S_{\tau_R\wedge \sigma_B})$.  Then 
\[
E_0 M_{\tau_R}=  \sum_{x\in B \; \mbox{{\footnotesize or}} \; |x|\geq R} p(x)a(x) +  \sum_{x\notin B, |x|<R} p(x) E_x M_{\tau_R} = \sum p(x)a(x) =1.
\]
Thus one has (\ref{a/Ea}).  For the proof of (\ref{eq_auB})  write  it as 
\beqn\label{a/a/a}
a^\dagger(x - \xi)=   E_x\big[a(S_{\tau_R}-\xi); \tau_R < \sigma_B\big] + E_x\big[a(S_{\sigma_B}-\xi); \tau_R \geq \sigma_B\big].
\eeqn
On passing to the limit as $R\to\infty$  the equality  (\ref{eq_auB}) then comes out in view of (\ref{u/a/Ha}), the last expectation converging to  $H^x_B\{a(\cdot-\xi)\} = a^\dagger(x-\xi) -u_B(x)$.

When  $X$ is   of finite range,  the identity (\ref{eq_auB}) (restricted to $x\notin B$)  is shown in the proof of     Proposition 4.6.3  of \cite{LL} (in a different way  from ours).
\end{remark}

\vskip1mm
Let $\hat H^x_B$ stand for the hitting distribution of $B$ for the dual (time-reversed) walk, in other words $\hat H^x_B(y) = H^{-x}_{-B}(-y)$ ($-B=\{-z: z\in B\}$).  Then  $H^x_B(x) =\hat H^x_B(x)$ and
\beqn\label{g/H}
g_B(x,y) = \hat H^y_B(x){\bf 1}_{\mathbb{Z}\setminus B}(y)  +\delta(x,y)  \quad \mbox{for}\;\; x\in B, y\in \mathbb{Z},
\eeqn 
where  ${\bf 1}_B$ is the indicator function of a set $B$.
\begin{lemma}\label{lem3.10}  Let  $\hat u_B$ be the dual of $u_B$: $\hat u_B(x) = a^\dagger(\xi-x) - \hat H^x_B\{a(\xi-\cdot)\}$ ($\xi\in B$). Then
\beqn\label{H/u/a}
H^x_B(y){\bf 1}_{\mathbb{Z}\setminus B}(x) =\hat u_B(y)+\sum_{z\in B} a(x-z) H^z_B(y)  - a^\dagger(x-y) \quad \mbox{for}\;\; x\in \mathbb{Z},  y\in B;
\eeqn
and  
\beqn\label{eqL3.10}
\begin{array}{rr}
{\displaystyle 1-H^x_B(\xi) = [1-H^\xi_B(\xi)]a^\dagger(x-\xi) +\sum_{z\in B\setminus\{\xi\}}\big[a(\xi-z)-a(x-z)\big]H^z_B(\xi)} \\
\mbox{for}\;\; \xi\in B,  x\notin B\setminus \{\xi\}.
\end{array}
\eeqn
\end{lemma}  
\begin{proof} Let  $x\in B$,  substitute the expression of  $g_B(x,y)$ given in  (\ref{g/H}) into (\ref{eqPS1}), rewrite the resulting identity in terms of the dual objects  and interchange $x$ and $y$, use the equality $\hat H^y_B(z) =H^z_B(y)$ for $y, z\in B$, and you obtain (\ref{H/u/a}).

For the proof of (\ref{eqL3.10})  we have only to  consider $x\neq \xi$, (\ref{eqL3.10}) being obviously true for $x=\xi$. Let  $\xi \in B$ and subtract  the equality  (\ref{H/u/a}) with  $x\notin B, y=\xi$ from that with  $x=y=\xi$.  Then one finds  
 \[
 -H^x_B(\xi) = - a(x-\xi) H^\xi_B(\xi) + \sum_{z\in B\setminus\{\xi\}} \big[a(\xi-z)-a(x-z)\big] H^z_B(\xi)   -1 +a^\dagger(x-\xi),
 \]
which  after a simple transposition of terms becomes (\ref{eqL3.10}).
\end{proof}

If $B$ is finite,  $\hat u_B(y) = \frac12  \lim_{x\to\infty}\big[H_B^{-x}(y)+ H_B^x(y)\big]$ and (\ref{H/u/a}) is given in \cite[P30.1]{S}.

The next lemma is  a consequence  of the second relation of Lemma \ref{lem3.10}.
\begin{lemma}\label{lem3.11}  Let $\sigma^2=\infty$.
If $B_n$, $n=1, 2,\ldots$ are  non-empty subsets of $\mathbb{Z}$ such that $\min \{|z|: z\in B_n\} \to \infty$ as $n\to\infty$ and    $P_0[\sigma_{-B_n}<\sigma_0]/P_0[\sigma_{B_n}<\sigma_0]$ is bounded, then
\beqn\label{Esc}
\lim_{n\to\infty} \frac{P_x[\sigma_{B_n}<\sigma_0]}{P_0[\sigma_{B_n}<\sigma_0]} = a^\dagger(x).
\eeqn
\end{lemma}
\begin{proof}  
If  $z, 0\in B$, then  $H^z_B(0) = H^0_{-B}(-z)$, and hence
$  \sum_{z\in B\setminus\{0\}}H^z_B(0)=1- H^0_{-B}(0) = P_0[\sigma_{-B} <\sigma_0].$
 Applying this as well as  identity (\ref{eqL3.10})  with  $\xi=0$, $B=B_n\cup \{0\}$ one deduces that  
\beqn\label{Bn/Bn}
\frac{P_x[\sigma_{B_n}<\sigma_0]}{P_0[\sigma_{B_n}<\sigma_0]}  =\frac{1-H^x_B(0)}{1-H^0_B(0)} = a^\dagger(x)  + \frac{\sum_{z\in B_n}[a(-z)-a(x-z)] H^z_{B_n\cup\{0\}}(0) }{P_0[\sigma_{B_n} <\sigma_0]}
\eeqn
of which the last ratio  tends to zero for each  $x$ under the condition imposed on $B_n$ in the lemma, for $a(-z)-a(x-z) \to 0$ ($-z\in B_n$) if $\sigma^2=\infty$. 
\end{proof}


\section{Proof of Theorem \ref{thm1.1} }

\begin{lemma}\label{lem4.1}\; 
Suppose $E Z=\infty$. Then\,  $a^\dagger(x)=H_{(-\infty,0]}^x\{a\}$.
\end{lemma}
\begin{proof}   We  consider only the case $x>0$,  the asserted formula for $x\leq 0$ being deduced from that for  $x>0$ by using (\ref{eq2.6}) as in (\ref{aa0}).

The  proof is based on  the fact that the function  $h(x):=a^\dagger(x)-H_{(-\infty,0]}^x\{a\}$ is non-negative and   harmonic on $x>0$ (according to Lemma \ref{lem3.9}). In view of the uniqueness of  harmonic function
it suffices to show
\beqn\label{eq4.1}
\liminf_{x\to\infty}\,\frac{a(x)}{f_r(x)} =0.
\eeqn

We have extended $f_r$  to a function on $\mathbb{Z}$, denoted also by $f_r$,  by (\ref{f_ext}), namely 
$f_r(x)= P[Z> -x]$ $(x\leq 0).$
Accordingly,  by  (\ref{f_Z})  we have 
$f_r(x)=\sum_{y=1}^\infty p(y-x)f_r(y)$ for all $x\in \mathbb{Z}$.

By the assumption of the lemma the walk is not right-continuous.  Hence  by Lemma \ref{lem3.1}
$\inf_{x<0} a(x)>0,$ so that for some  constant $C$
\beqn\label{eq4.2}
 f_r(x)\leq Ca(x) \quad  \mbox{for}  \;\; x<0.
 \eeqn

Define the operators $P$ and $P^-$ by
\[
Pf(x)=\sum_{y\in\mathbb{Z}} p(y-x)f(y) \quad \mbox{and} \quad P^- f(x)=\sum_{y\leq 0} p(y-x)f(y)\qquad (x\in \mathbb{Z}), 
\]
respectively.  Put $G_n(x,y)=p^n(y-x)+\cdots+p(y-x)+\de(x,y)$ and let $G_n$ also denote the corresponding operator. We may suppose $\inf_{x>0} a(x)>0$, otherwise $a(x)$ vanishing for all $x>0$ so  that (\ref{eq4.1}) is plainly evident.
Owing to    (\ref{eq4.2})   relation (\ref{eq4.1}) then  follows if we can show
\beqn\label{eq4.3}
\lim_{n\to\infty}\, \frac{P^n a^\dagger(0)}{P^n f_r(0)}=0,
\eeqn
for if (\ref{eq4.1}) does not hold,  $a^\dagger =a+\de(\cdot,0)$ must dominate a positive multiple of  $f_r$ so that  (\ref{eq4.3}) is impossible.   

From the identity $Pa=a+\de(\cdot,0)$ one deduces   by induction  that
\beqn\label{W1}
P^na^\dagger(x)= a(x)+G_n(x,0).
\eeqn
On the other hand one obtains  that   $Pf_r= f_r+ P^-f_r$ and  by induction again  $P^nf_r(x)= f_r(x)+ G_{n-1} P^- f_r(x)$, which can be rewritten as
\beqn\label{W2}
P^nf_r(x)= f_r(x){\bf 1}_{[1,\infty)}(x)+  \sum_{y\leq 0}G_{n}(x,y)f_r(y).
\eeqn
Since 
$ \sum_{y\leq 0} f_r(y) =
 \sum_{y\geq 0} P[Z>y]=\infty$ due to  the assumption of the  lemma,  for any $K>0$ one  can  choose  a positive integer $M$ so that 
$\sum_{-M\leq y\leq 0} f_r(y)\geq K;$ 
and hence
\[
P^nf_r(0)\geq K\min_{0\leq z \leq M} G_{n}(0,-z)\geq 2^{-1}K G_{n}(0,0)
\]
if $n$ is large enough,  for the recurrence of the walk  implies $\lim_{n\to\infty}G_n(z,0)/G_n(0,0)=1$ (cf.\cite[P2.6]{S}).  Combined with (\ref{W1})    the  inequality   derived above    implies that
$P^n f_r(0)\geq \frac12 K P^{n} a^\dagger(0)$ for all sufficiently large $n$ and     we can conclude the required relation (\ref{eq4.3}). 
\end{proof}


\textit{ Proof of Theorem \ref{thm1.1}}.  By Lemma  \ref{lem3.9}  the formula of (i)  follows if we verify  its  special case $y=0$. Note that  $g_{(-\infty,0]}(x,0) =\delta(x,0)$. It is then obvious that  if $EZ=\infty$,  (i) and  (iii)  follow from  Lemma \ref{lem4.1} and  (\ref{eq4.1}), respectively.  
 Let $EZ<\infty$. Then by Lemma  \ref{lem3.9}
 the difference  $a^\dagger(x)- H_{(-\infty,0]}^x\{a\}$ is non-negative and harmonic on $x>0$, so that it is a constant multiple of $f_r(x)$.     The constant factor is determined by using 
 Lemmas  \ref{lem3.5},  \ref{lem3.7}, and \ref{lem3.8}. (Note that $\bar a(x) \sim a(x)$ if $\sigma^2<\infty$.)  By these  the second assertion (ii)  also follows. 
 \qed

\section{ Proof of Theorem \ref{thm1.2} }

\quad \textit{Proof of} (i).   The first half of   (i) of Theorem \ref{thm1.2} follows from Lemma \ref{lem3.6} and Corollary  \ref{cor1}. 
Suppose that $ EZ<\infty$ and  $\sigma^2=\infty$. 
 The formula (\ref{double}), what is asserted in the second half, may be written as
\beqn\label{double4}
{a(-x)}\,  \sim \,  \frac1{E Z} \, \sum_{k=1}^\infty P\big[k<Z\leq x+k\big]a(k) \quad (x\to\infty).
\eeqn
The hitting distribution of the half line $[0,\infty)$ by the random walk started at a negative 
site $-x$ agrees with that by the ascending ladder height process  started at $-x$.
Since    $G(x',x''):=u_{\rm as}(x''-x')$ ($x'\leq  x''$) is the Green function of the ladder height process  it accordingly  follows that
 \beqn\label{rep_H}
 H^{-x}_{[0,\infty)}(k) = \sum_{y=1}^xu_{\rm as}(x-y)P[Z=k+y] = \sum_{w=0}^{x-1}u_{\rm as}(w)P[Z=k+x-w]. 
 \eeqn 
Since for  each  $w$,  $\sum_{k=0}^\infty P[Z=k+x-w] a(k) 
 \to 0$ as $x\to\infty$ and $u_{\rm as}(w)\to 1/EZ$ ($w\to\infty$),  we can conclude  (\ref{double4}) owing to 
  Corollary  \ref{cor1}(ii) that gives the identity $H_{[0,\infty)}^{-x}\{a\}= a(-x)$.

\vskip2mm
\textit{ Proof of} (ii).\, If  $ \limsup_{x\to-\infty}a(x)<\infty$, then $H_{(-\infty,0]}^x\{a\}$ is bounded, so that $EZ$ cannot be finite, for otherwise by Corollary  \ref{cor1}(ii)   $a(x)+a(-x)$ must be bounded  which  is impossible  in view of  (\ref{a1}). This shows (ii) by virtue  of Lemma \ref{lem3.1}.

\vskip2mm
\textit{ Proof of}  (iii).\;
Let $\sigma^2=\infty$ and $E Z<\infty$.  Then $E|\hat Z|=\infty$ so that $a(-x)=H_{[0,\infty)}^{-x}\{a\}$ according to Corollary \ref{cor1} and  by the first equality of (\ref{rep_H})  we have
\beqn\label{a/b}
a(-x) =\sum_{y=1}^x u_{\rm as}(x-y)b(y),
\eeqn
where 
\beqn\label{D_b}
 b(y) = \sum_{k=1}^\infty P[Z= y+k]a(k).
 \eeqn
Since $b(y)\to 0$ under $EZ<\infty$, (\ref{a/b}) yields
\beqn\label{double02}
a(-x)\sim 
 \frac1{E Z}\sum_{y=1}^x b(y).
\eeqn
  By (\ref{sharp})
\beqn\label{Epr_b}
b(y) = \sum_{k=0}^\infty\sum_{w=0}^\infty v_{\rm ds}(w)p(y+k+w)a(k),
\eeqn
which by a change of variables one can rewrite as
\[
b(y)=\sum_{j=1}^\infty p(y+j) 
\sum_{w= 0}^{j}v_{\rm ds}(w)a(j-w).
\]
 Since  $f_r$  is sub-additive (see (\ref{sb-add})) and  $a(x)EZ \sim f_r(x)=  v_{\rm ds}(0)+\cdots+v_{\rm ds}(x-1)$,  the inner sum is   dominated from below and above by  positive multiples of $[f_r(j)]^2$,
 so that
$b(y)\asymp  \sum_{j=1}^\infty p(y+j) [f_r(j)]^2$.
 Now substitution into (\ref{double02})  yields
\[
a(-x) \asymp \sum_{y=1}^x  \sum_{j=1}^\infty p(y+j) [f_r(j)]^2.
\]
 In view of (\ref{c<c}) (or Lemma \ref{lem3.6}(ii))  one can replace $f_r(j)$  by $j/m_-(j)$, showing  (\ref{a00}), the desired asymptotics of $a(-x)$.    The rest of (iii) is readily ascertained to be true by  (\ref{a00}), Lemma \ref{lem4.1} and (\ref{H<H}). \qed

\begin{lemma}\label{lem5.1} \,   If $EZ =\infty$ and  $a(-x)$ is almost increasing in the sense that $a(-y)\geq \de a(-x)$ if $y> x \geq x_0$ with some $\de>0$
and $x_0>0$, then 
\[
(*)\qquad a(x)/f_r(x)\to 0 \quad \mbox{as}\quad x\to\infty. \qquad\qquad
\]
\end{lemma}
\begin{proof}   If  $E|\hat Z|<\infty$, then $f_r(x)\sim Cx$ and the assertion of the lemma is obvious. Let $E|\hat Z|=EZ=\infty$
and we show that if $a(x)$ is almost increasing, then $a(-x)/f_l(x) \to 0$,  which  by duality  amounts to the same as the assertion of the lemma.  
For $x\geq 1$, $a(-x) =H^{-x}_{(\infty,0]}\{a\}$ so that we have (\ref{a/b}) and for the present purpose 
it suffices to show that $b(y)\to 0$. 
 Rewrite the expression of $b(y)$ in (\ref{Epr_b}) as
 \[
 b(y) = \sum_{w=0}^\infty v_{\rm ds}(w) \sum_{k=0}^\infty a(k)p(y+k+w).
 \]
 Now suppose that $a(x)$ 
is almost increasing.  Then
 \beqn
  \sum_{k=0}^\infty a(k)p(y+k+w) \leq C \sum_{k=y}^\infty a(k)p(k+w) + \sum_{k=0}^{x_1}a(k)p(y+k+w)
  \eeqn
   for some $x_1$. On noting  $ \sum_{w=0}^\infty v_{\rm ds}(w) \sum_{k=0}^\infty a(k)p(k+w)= b(0)= Ea(Z) <\infty$ according to (\ref{D_b}),  the dominated convergence therefore concludes  
that $b(y)\to 0$,  as desired. 
\end{proof}

\begin{remark}\label{rem5.1}  We apply Lemma \ref{lem5.1} 
 to verify the last assertion of Remark \ref{rem2.1}(e).
 Let $X$ belong to the domain of attraction of a Cauchy distribution (a stable law with exponent $\alpha=1$ and  skewness parameter $\beta=0$).  Then $\bar a(x)$ is dominated by a slowly varying function  \cite[Remark 62(ii)]{Upot}.  Suppose $\lim  P[S_n>0]=\rho$  in addition.  Then  $a(-x)/a(x) \to1$ and  $f_r(x)$ is regularly varying with index  $1-\rho$ \cite{R}, so that  it  plainly follows  that $\bar a(x)/f_r(x)\to 0$  if $\rho<1$. Let  $\rho =1$. Then $a(x)\sim  a(-x) \sim   \int_{x_0}^x  F(-s)\{A^2(s)\}^{-1}ds$ for some $x_0>0$, where $A(s)=\int_0^s\big[1-F(t)-F(-t)\big]dt$  (cf. \cite[Theorem 7]{Upot}), in particular $a(-x)$ is almost increasing.   Thus Lemma \ref{lem5.1} verifies $(*)$, hence $\bar a(x)/f_r(x)\to 0$.     We also know that  $f_r(x) \sim 1/\int_x^\infty \big[F(-s)\big/\ell^*(s)\big]ds$, where $\ell^*(s)=\int_0^s P[Z>t]dt$ \cite[Lemma 3.1]{Uexit}). However, it is not clear whether   $(*)$ can be  deduced  directly from  these asymptotic relations. 
\end{remark}


\section{Proofs of Propositions \ref{prop1} and \ref{prop2}} 
 
We employ the following identities:
\beqn\label{eq4.4}
P_x[\sigma_R<T]=\frac{g_{(-\infty,0]}(x,R)}{g_{(-\infty,0]}(R,R)};\quad P_x[\sigma_R<\sigma_0]=\frac{a^\dagger(x)+a(-R)-a(x-R)}{2\bar a(R)}.
\eeqn
    If either $a(-x)/a(x)\to 0$ or $a(x)/a(-x)\to 0$  ($x\to \infty$), then  by Lemma \ref{lem3.5}
  $g_{(-\infty, 0]}(R,R) \sim  2\bar a(R)$
so that
   \beqn\label{C_L6.1}
   \frac{P_x[\sigma_R<T\,]}{P_x[\sigma_R<\sigma_0]} 
   \sim \frac{g_{(-\infty,0]}(x,R)}{g_{\{0\}}(x,R)}. 
   \eeqn
According to \cite[Theorem 4]{Upot}  $a(\mp x)/a(\pm x)\to 0$ if $m_\pm(x)/m(x) \to 0$. 

\begin{lemma}\label{lem6.1}  If $a(-x)/a(x) \to 0$  ($x\to\infty$),  then  for each $M\geq 1$, 
   uniformly for $-M < z < y$
\[
a(-y)-a(z-y)= o(a(z)\vee1)\qquad  (y\to\infty);
\]
 in particular  as $R\to\infty$
\[
 P_x[\sigma_R<\sigma_0] \,\sim  \,a^\dagger(x)/a(R) \quad \mbox{uniformly for \;\; $-M<x<R$}.
 \]
 \end{lemma}
\begin{proof}   This---verified readily by using Lemma \ref{lem3.2}---is contained in Lemma 37(ii) of \cite{Upot}. 
\end{proof}
 
{\bf 6.1} {\sc Proof of Proposition \ref{prop1}.}  The first case of (\ref{E2}) follows from the second equivalence in (\ref{E11}) together with Lemma \ref{lem6.1}.   Without recourse to  (\ref{E11}) it may be verified as follows.
 If   $E Z<\infty$,  by the expression of $g_{(-\infty,0]}(x,y)$ in Theorem A(ii)  it follows  that  uniformly for $x<R$ as $R\to\infty$
\beqn\label{g1}
g_{(-\infty,0]}(x,R) \sim  f_r(x)/E Z
   \eeqn
 and the relation (\ref{E2}) asserted  in  Proposition \ref{prop1}  follows from  (\ref{C_L6.1}) and Lemma \ref{lem6.1}; 
  the uniformity of the convergence is assured by  $\lim f_r(x)/a(x) =EZ$.   
  
  The case $EZ=\infty$ of (\ref{E2})  is essentially contained in   the next lemma.
\begin{lemma}\label{lem6.2} \, If $\sigma^2=\infty$, 
 ${\displaystyle \lim_{R\to\infty} 2\bar a(R)P_x[T<\sigma_R<\sigma_0]= H_{(-\infty,0]}^x\{a\} \quad ( x\in \mathbb{Z}).}$
\end{lemma}
\begin{proof}  First of all it is pointed out that we may suppose $P[X\geq 2]>0$ so that for all $x$, 
so that $P_x[\sigma_R<\sigma_0]>0$, for otherwise we have $H^x_{(-\infty,0]}\{a\} =0$ for all $x$ and  the result is obvious. 
Now  observe that  the equality  
$ P_x[T<\sigma_R<\sigma_0]= P_x[\sigma_R<\sigma_0] - P_x[\sigma_R< T \,]$  together with  (\ref{eq4.4})  
 yields
\beqn\label{eq4.6}
2\bar a(R)P_x[T<\sigma_R <\sigma_0] =\big[a^\dagger(x)+a(-R)-a(x-R)\big]\bigg(1-\frac{P_x[\sigma_R<T\,]}{P_x[\sigma_R<\sigma_0]} \bigg) 
\eeqn
on the one hand,  and  since $P_x[T<\sigma_R<\sigma_0]=\sum_{y<0}P_x[S_T=y, T<\sigma_R]P_y[\sigma_R<\sigma_0]$,
\beqn\label{eq4.7}
2\bar a(R)P_x[T<\sigma_R<\sigma_0]=\sum_{y<0}P_x[S_T=y, T<\sigma_R]\big[a(y)+a(-R)-a(y-R)\big]
\eeqn
on the other hand. 

Suppose $E Z<\infty$  so that the first case of  Proposition \ref{prop1} is applicable. Then under $\sigma^2=\infty$    
  the RHS of  (\ref{eq4.6}) converges to
$a^\dagger(x) - Af_r(x)/EZ$, hence   the identity of the lemma in view of Corollary \ref{cor1}(i).

  By  letting $R\to\infty$ in the  identities (\ref{eq4.6}) and (\ref{eq4.7}), with the help of Fatou's lemma for the infinite series on the RHS of (\ref{eq4.7}),    we obtain
\begin{eqnarray} \label{eq4.71}
H_{(-\infty,0]}^x\{a\}&\leq& \liminf_{R\to\infty}  2\bar a(R)P_x[T<\sigma_R<\sigma_0] \nonumber\\
&=& a^\dagger(x)-a^\dagger(x)\limsup_{R\to\infty} \frac{P_x[\sigma_R<T\,]}{P_x[\sigma_R<\sigma_0]} \leq a^\dagger(x),
\end{eqnarray}
provided $\sigma^2=\infty$.  If $E Z=\infty$, then
  the two extreme  members  in (\ref{eq4.71}) must coincide  owing to Corollary \ref{cor1}(ii),
entailing that the two inequalities above are the equality, of which the latter means that the $\limsup$ vanishes---showing the relation   of  Proposition \ref{prop1}.  We can interchange the $\liminf$ 
and the  $\limsup$ in (\ref{eq4.71}),  which gives  the equality  of the lemma.
  \end{proof}
   
   \begin{remark}\label{rem6.1} If  $\sigma^2<\infty$,  we have 
   $
{P_x[\sigma_R<T\,]}/{P_x[\sigma_R<\sigma_0]}  \sim {g_{(-\infty,0]}(x,R)}/{g_{\{0\}}(x,R)}$ uniformly for $x\in \mathbb{Z}$.
 By (\ref{eq4.6}) and (\ref{g1}) it therefore follows that as $R\to\infty$
 \[
 2\bar a(R)P_x[T<\sigma_R<\sigma_0]= g_{\{0\}}(x,R) - \frac{f_r(x)}{EZ} \{1+o(1)\} \, \longrightarrow\, a^\dagger(x) +\frac{x}{\sigma^2} - \frac{f_r(x)}{EZ}
 \]  
 ($ x\in \mathbb{Z}$), where the equality is uniform for $x<R$ (but the convergence is not).
  \end{remark}
  
 {\bf 6.2} {\sc Proof of Proposition \ref{prop2}.} We prove the dual assertion that follows.
   
\begin{lemma}\label{lem6.3}\,  If  $E |\hat Z|< \infty$,  then, uniformly for $x\in \mathbb{Z}$, 
as $R\to\infty$
\[P_x\big[\sigma_{[R,\infty)}< T \big] = P_x\big[\sigma_{R}< T \big] \{1+o(1)\}
\]
and  as $x\wedge R \to\infty$ under  $ x < R$
\beqn\label{E1}
P_x\big[\sigma_{[R,\infty)}< T \big] \sim  \frac{f_l(R)- f_l(R-x)}{f_l(R)}\leq \frac{f_l(x)}{f_l(R)}.
\eeqn
\end{lemma}

\begin{proof}   Let $\sigma^2=\infty$. Then the summability   of  $\hat Z$ implies $a(R)/a(-R)\to0$  as well as  the tightness of  the family  $\{H_{(-\infty,0]}^y: y>0\}$,  which together  imply  that
for each $z$ the probability
\[
P_{R-z}[T<\sigma_R] =P_{R-z}[\sigma_0<\sigma_R] + P_{R-z}[T<\sigma_R <\sigma_0]
\]
tends to zero (with  $z$ fixed)---use  (\ref{eq4.4}) for the first term on the RHS; note that the second term is less than  $\sum_{y<0} H^{R-z}_{(-\infty,0])}(y)P_y[\sigma_R<\sigma_0]$. Hence
\[
\sup_{y>R}P_y\big[T<\sigma_R\big] =\sup_{y'>0}\sum_{z > 0} H_{(-\infty,0]}^{y'}(-z)P_{R-z} \big[ \tilde T <\sigma_R\big]\longrightarrow 0 \quad\mbox{as} \;\; R\to\infty,
\]
where  $\tilde T= \inf \{n\geq 0: S_n\leq 0\}$, and hence the ratio
\[ 
\frac{P_x[\sigma_{[R,\infty)}< T] - P_x[\sigma_R<T]}{P_x[\sigma_{[R,\infty)}< T] } 
=\sum_{y>R} P_x\big[S_{\sigma_{[R\infty)}}=y\,\big|\, \sigma_{[R,\infty)}<T \big]P_y\big[T<\sigma_R\big] 
\] 
tends to zero  uniformly for $x>0$, which shows the first half of the lemma. 

For the proof of  the second half
we derive the asymptotic form of $P_x[\sigma_R<T]$  by using (\ref{eq4.4}).   Let   $E |\hat Z| < \infty$. Then   one   obtains  $g_{(-\infty,0]}(R,R) \sim f_l(R)/E|\hat Z|$ because  of  the dual of (\ref{g1}),  so that
\[
P_x\big[\sigma_{[R,\infty)}< T \big] \sim P_x\big[\sigma_R<T \big]  = \frac{\sum_{k=0}^{x-1} v_{\rm ds}(k)u_{\rm as}(R-x+k)}{f_l(R)}\big(E|\hat Z| +o(1)\big).
\]
In order to verify (\ref{E1})   it suffices to see that for each $K>0$, $\sum_{k=0}^K u_{\rm as}(y+k)$ divided by  $\sum_{k=0}^{x-1} u_{\rm as}(y+k)$ tends to zero as $ x\to\infty$ uniformly for $y\geq0$.
By  Lemma \ref{lem6.4} below  it follows that for each  $k\leq K$ and $j\geq k$,   $u_{\rm as}(y+j) \geq u_{\rm as}(y+k)u_{\rm as}(j-k)$ so that  
\[\sum_{j=0}^{x-1} u_{\rm as}(y+j)\geq  u_{\rm as}(y+k)[u_{\rm as}(0)+\cdots +u_{\rm as}(x-k-1)]
\]
 and hence  the ratio in question is dominated by $K/ f_l(x-K)$, which tends to zero as required. The inequality in (\ref{E1}) follows by the  sub-additivity of $f_l$   (cf. (\ref{sb-add})). 
\end{proof}

\begin{lemma}\label{lem6.4} \, For all integers $x, y \geq 0$,\;
$
u_{\rm as}(x+y)\geq    u_{\rm as}(x)u_{\rm as}(y).
$
\end{lemma}
\begin{proof}   The ratio $u_{\rm as}(x+y)/u_{\rm as}(x)$ is not less than the conditional probability that $x+y$ is an ascending ladder point given so is $x$, but this conditional probability equals $u_{\rm as}(y)$, showing the inequality of the lemma. 
\end{proof}
 
\begin{remark}\label{rem6.2}
Let $E|\hat Z|<\infty$ and $\sigma^2=\infty$.  Then by Corollary \ref{cor1} $\bar a(R)  \sim  f_l(R)/E|\hat Z|$, 
so that by  (\ref{C_L6.1}) and Lemma \ref{lem6.3} it holds that whenever  $x<R$ and $x\to \infty$ 
\[
\frac{P_x [\sigma_{R}<  T ] }{P_x\big[\sigma_{R}< \sigma_{0} \big]} \sim \frac{[f_l(R)-f_l(R-x)]/E|\hat Z|}{a(x)+a(-R)-a(x-R)}. 
\]
For any $0<\ep<1/2$,  the RHS can be shown to approach unity uniformly for $\ep R<x<R$   if it is further supposed that  $1-F(x)$ is regularly varying with index $ -\alpha < -1$ and $F(-x)/[1-F(x)] \to 0$
 (cf. \cite[Section 8.1.1]{Upot})---if $\alpha =1$,  it may  approach  zero uniformly for
 $0<x< (1-\ep_R)R$ for an appropriate $F$ and $\ep_R>0$ decreasing to 0 [this actually takes place  if, e.g., $P[X= x] = x^{-2}(\log x)^{-\la}\{1+ O(x^{-1})\}$,  
 $F(-x) = x^{-1}(\log x)^{-\la-\de}\{1+O(x^{-1})\}$  (as $x\to\infty$) and $\ep_R \sim \exp\big\{-(\log R)^{\ep}\big\}$ 
 with $\la >1$, $0<\de<1$ and $0<\ep <\de \wedge \frac12$---proof is quite involved even in such a  particular  case]. 
 \end{remark}

\section{The random walk conditioned on $\sigma_0=\infty$}

Write $\tilde P_x$ for $P_x[\,\cdot\, |\sigma_0=\infty]$~($x\in \mathbb{Z}$), the probability law of the  conditional process $S_n$ given that it never
visits the origin.  It is defined as a  limit law of $P_x[\cdot |\sigma_0>k] $ as $k\to\infty$. If $\sigma^2=\infty$,  suppose  $P[X\leq -2]P[X\geq 2] > 0$ so that $a^\dagger(x)>0$ for all $x$. 
[If $P[X\leq -2]P[X\geq 2] =0$, this conditioning   forces    the walk  to stay either the positive or negative half line once it get into there, yielding the process  represented by the harmonic transform by means of  $f_r$ or $f_l$   according as the starting site is positive  or negative, respectively.] 
The conditional process is Markovian with state space $\mathbb{Z}\setminus\{0\}$ and the $n$-step transition law   given by
\beqn\label{eq5.1}
\frac1{a(x)}q^n(x,y)a(y) \;\; (x, y\neq 0) \quad \mbox{where} \;\; q^n(x,y)=P_x[S_n=y,  n<\sigma_0].
\eeqn
   Indeed, for $n<k$   
\beqn\label{c_trans_l}
P_x[ \,S_n=y\, |\, \sigma_0>k] =q^n(x,y)\frac{P_y[\sigma_0>k-n]}{P_x[\sigma_0>k]},
\eeqn
and as $k\to\infty$,  $P_y[\sigma_0>k]/P_x[\sigma_0>k]\to a(y)/a(x)$ while ${P_y[\sigma_0=k-n]}/{P_x[\sigma_0>k]} \to 0$ for each $n$  (see \cite[T32.1, T32.2]{S}) so that the ratio in the RHS of (\ref{c_trans_l}) converges to $a(y)/a(x)$, showing (\ref{eq5.1}). Let $B(R) = \mathbb{Z}\setminus (-R,R)$. Then by Lemma \ref{lem3.11} 
 \[
 \lim_{R\to\infty} \frac{P_x[\sigma_{B(R)}<\sigma_0]}{P_0[\sigma_{B(R)}<\sigma_0]} \,  = a^\dagger(x)
 \]
(in case $\sigma^2< \infty$  see (\ref{Bn/Bn}) or (\ref{eq5.7}) below),  from which one can easily deduce that
  \beqn \label{S7}
\mbox{the conditional law }\; P_x[\, \cdot  \, | \,  \sigma_{B(R)}<\sigma_0] \;\;\mbox{converges to} \;\;  \tilde P_x\;  \;  \mbox{as}\;\; R \to\infty
 \eeqn
in the sense of  convergence of   finite dimensional distributions.

It  follows from (\ref{eq5.1})  that
\beqn\label{eq5.2}
\tilde H_{(-\infty,0]}^x(y):=\tilde P_x[\,S_T=y]=\frac1{a^\dagger(x)}H_{(-\infty,0]}^x(y)a(y)\quad\quad (x\in \mathbb{Z}, y<0).
\eeqn
Therefore by Corollary  \ref{cor1} 
\beqn\label{eq5.3}\tilde P_x[\,T<\infty] = 1-\frac{Af_r(x)}{E[Z] a^\dagger(x)} \quad  \quad (x\in \mathbb{Z}).
\eeqn
[Recall  $f_r(x) = P[Z> -x]$ for $x\leq0$.]  (\ref{eq5.1}) also shows $\sum_n \tilde P_x[S_n =y] <\infty$.   Hence
 \beqn\label{eq5.4}
\tilde P_x[\,|S_n|\to\infty \;\;\mbox{as} \;\; n\to\infty]=1.
\eeqn
In fact we have that in case  $\sigma^2=\infty$, for every $x\in \mathbb{Z}$ 
    \beqn\label{eq5.5}
\begin{array}{ll}
{\rm(a)}\quad 
\tilde P_x[\,\lim S_n= +\infty]=1 \quad &\mbox{if \,\; $ E Z<\infty$},\\[1mm]
{\rm(b)}\quad 
\tilde P_x[\,\limsup S_n= +\infty\, \mbox{and} \,\liminf S_n=-\infty]=1 \;\; &\mbox{if\,\,$ E Z= -E \hat Z=\infty$;}
\end{array}
\eeqn
and in case $\sigma^2<\infty$,  either $\lim S_n= +\infty$ or $\lim S_n  =-\infty$ with $\tilde P_x$-probability one and
\beqn\label{eq5.6}
\tilde P_x[\,\lim S_n= + \infty]=\frac{a^\dagger(x)+\sigma^{-2}x}{2a^\dagger(x)} \quad\quad (x\in \mathbb{Z}).
\eeqn

The  two identities in (\ref{eq5.5})  are readily deduced from (\ref{eq5.3}) (or rather (\ref{eq5.2})) and  its dual relation as well as  (\ref{eq5.4}) by virtue of  Corollary \ref{cor1}.   
It is also noted that  for each $M>1$
\[
\begin{array}{ll}
 \tilde P_x[\sigma_R< \infty ] \to 1\;\;\;  \mbox{as $R\to\infty$  uniformly for $x\in (-M,R)$}\;\; \mbox{ if}\;\; a(-x)/a(x)\to 0; \mbox{and}\\[1mm]
\tilde P_x[\sigma_R< T] \to \left\{\begin{array}{ll} \! 1\;\;\;\;  \mbox {as $x\to \infty$ uniformly for $R>x$}\;\; &\mbox{if}\;\;  EZ<\infty,\\
\! 0 \;\;\;\; \mbox {as $R\to\infty$ for each  $x\in \mathbb{Z}$} \;\; &\mbox{if}\;\;  EZ=\infty,
\end{array}\right.
\end{array}
\]
which together in particular show (a). Here the first relation follows from Lemma \ref{lem6.1}  in view of  $\tilde P_x[\sigma_R<\infty] =P_x[\sigma_R<\sigma_0]a(R)/a^\dagger(x) $ and  the second from Proposition \ref{prop1}, which shows that if $EZ<\infty$,  then uniformly for $0\leq x<R$, 
\[
\tilde P_x[\sigma_R< T] =P_x[\sigma_R< T]a(R)/a^\dagger(x) \sim Af_r(x)/ [a^\dagger(x) EZ] \quad (R\to\infty).
\] 

The formula (\ref{eq5.6})  is obtained by applying  a theorem from the theory of Martin boundary  (see \cite[Theorem III29.2]{RW}: the Martin kernel $\kappa(\cdot,\pm)$ relative to a reference point $\xi\in \mathbb{Z}\setminus\{0\}$ is given by $[a(\cdot) \pm \sigma^{-2} \,\cdot\,]/ [a(\xi)+ \sigma^{-2}\xi]$). The conditional process $(\tilde P_x)_{x\neq 0}$  is a harmonic transform of the walk with absorption  at the origin whose Martin  boundary  contains  exactly two extremal  harmonic functions $h_+$ and $h_-$ given by  $h_{\pm}(x)=\lim_{y\to\pm \infty} g_{\{0\}}(x,y)/\sum_{z\neq 0} p(z)g_{\{0\}}(z,y)
=a(x)\pm\sigma^{-2}x $ ($x\neq 0$), provided $\sigma^2<\infty$. It is noticed that if $\sigma^2=\infty$, there is only one harmonic function, hence a unique Martin boundary point: $\lim_{|y|\to\infty}g_{\{0\}}(\cdot, y)/g_{\{0\}}(\cdot,\xi)  =a(\cdot)/a(\xi)$, so that two geometric boundary points  $+\infty$ and $-\infty$ are not distinguished in the Martin boundary 
whereas the walk itself discerns them provided that either $E Z$ or $E \hat Z$ is finite.   

The RHS of  (\ref{eq5.6})  equals the limit as $R\to\infty$  of 
$P_x[ S_{\sigma [R,\infty)} >0\,|\, \sigma_{B(R)}<\sigma_0]$ 
(the probability of the walk exiting the interval $(-R,R)$ from the upper boundary)
---as is shown by (\ref{eq5.7}) below, and prompted
by this fact  we here provide a direct proof of  (\ref{eq5.6}) that is based on (\ref{S7}).  

Suppose $\sigma^2 <\infty$. 
Using (\ref{a3})  one  infers  first 
$P_x[\sigma_{[R,\infty)}\vee \sigma_{(-\infty,-R]}<\sigma_0] =o(1/R)$ 
and then as $R\to\infty$,  
\beqn\label{eq5.7}
P_x[\sigma_{[R,\infty)} <  \sigma_0] \,\sim\,\frac{a^{\dagger}(x)+\sigma^{-2}x}{2\bar a(R)},\quad 
P_x[\sigma_{B(R)}<\sigma_0]\,\sim\,\frac{a^{\dagger}(x)}{\bar a(R)}.
\eeqn

 By  (\ref{eq5.4}) it suffices to show
\beqn\label{f_EQ} 
\lim_{M\to\infty} \tilde P_x[\sigma_{[M,\infty)} <\sigma_{(-\infty,-M]}]
=
\lim_{R\to\infty} P_x\big[S_{\sigma [R,\infty)} >0\,\big|\, \sigma_{B(R)}<\sigma_0\big],
\eeqn
the LHS being equal to $\tilde P_x[S_n\to\infty]$.
For verification of  (\ref{f_EQ})  put   $\tau^-_r= \sigma_{(-\infty,-r]}$ and  $\tau^+_r= \sigma_{[r,\infty)}$ for $r>1$ and denote by 
${\cal E}_R $ the event $\{\sigma_{B(R)}<\sigma_0\}$. 
Then one can readily see that
\beqn\label{e11}
\limsup_{R\to\infty}P_x[\tau^+_M<\tau^-_M, \tau^-_R<\tau^+_R\,|\, {\cal E}_R] \leq P_x[\tau^+_M<\tau^-_M <\infty\,|\, {\cal E}_R],
\eeqn
of which the RHS obviously tends to zero as $M\to\infty$.
Similarly
\beqn\label{e12} \limsup_{R\to\infty}P_x[\tau^-_M<\tau^+_M, \tau^+_R<\tau^-_R\,|\, {\cal E}_R] \, \longrightarrow\, 0 \quad (M\to\infty).
\eeqn
Thus,  for  $M$ and $R =R(M)$ chosen large enough, 
$ P_{x}[\tau^+_M<\tau^-_M, \tau^{+}_R<\tau^{-}_{R}\,|\, \sigma_{B(R)}<\sigma_{0}]$
approximates closely to  $P_{x}[\tau^{+}_M<\tau^{-}_M \,|\,{\cal E}_R] $  because of (\ref{e11}) and to 
$P_{x}[ \tau^{+}_R<\tau^{-}_{R}\,|\, {\cal E}_R]$ because of (\ref{e12}). It therefore follows that 
$$\lim_{R\to\infty} \big| P_{x}[\tau^{+}_M<\tau^{-}_M \,|\,{\cal E}_R] - P_{x}[ \tau^{+}_R<\tau^{-}_{R}\,|\, {\cal E}_R]\big| 
\, \longrightarrow\, 0  \quad\mbox{as}\;\;\; M\to\infty,$$
 which verifies  (\ref{f_EQ}) since  the above limit equals $\big| \tilde P_{x}[\tau^{+}_M<\tau^{-}_M \,] -  {\displaystyle \lim_{R\to\infty}  P_{x}[ \tau^{+}_R<\tau^{-}_{R} ]\,|\, {\cal E}_R] \big| }$. \qed

\section{Appendix}
 Put $Z'= S_{\sigma[S_0,\infty)} -S_0$, the weak ascending  ladder height. The renewal functions for the strictly and weakly ascending ladder height processes are defined by
$U_{{\rm as}}(x) =1+ \sum_{k=1}^\infty P[ Z_1+\cdots +Z_k\leq x]$  and 
$ V_{{\rm as}}(x) = 1+ \sum_{k=1}^\infty P[ Z'_1+\cdots +Z'_k\leq x]$ ($x=0, 1, 2, \ldots$).
Here  $(Z_n)$ and $(Z_n')$ are  i.i.d.\,copies of  $Z$ and $Z'$, respectively.  It follows  \cite[Section XII.1]{F} that $P[Z' \leq x] = P[Z' =0] + P[Z'>0] P[Z\leq x]$ and
\[  V_{{\rm as}}(x) =  U_{{\rm as}}(x)/P[Z'>0].
\]

Let $\tau =\sigma_{[1,\infty)}$,  $\tau' =\sigma_{[0,\infty)}$ and $c(t) = e^{-\sum_{1}^\infty k^{-1}t^kp^k(0)}$ ($t \geq 0$). Then $S_{\tau'} \stackrel{{\rm law}}=Z'$ and $S_\tau \stackrel{{\rm law}}= Z$ under $P_0$ and 
 $1- E_0[t^{\tau'}z^{S_{\tau'}}] = c(t)(1- E_0[t^{\tau}z^{S{\tau}}])  $ for $0\leq t <1, 0< |z|<1$ (\cite[Proposition 17.5]{S}, \cite[Section XVIII.3]{F}),
so that on  letting $z\downarrow 0$ and $t\uparrow 1$ in this order
\[
P[ Z' >0] =1/V_{\rm as}(0) =c(1) =c.
\]

For $x=1,2,\ldots,$ put   $\tau(x) = \inf\{n\geq 1: Z_1+ \cdots + Z_n \geq  x\}$,  the first epoch when the ladder height process enters  $[x,\infty)$. Then  $P[\tau(x)>n]=P[Z_1+\cdots+Z_n\leq x-1]$ ($n\geq 1$)  and hence
\beqn\label{sb-add}
U_{{\rm as}}(x-1) = E\tau(x) \quad (x=1,2, \ldots),
\eeqn
which especially shows that $f_l(x)$, which equals $c \, U_{{\rm as}}(x-1)$, $x\geq 1$,   is sub-additive.

\end{document}